\documentclass[12pt,oneside,reqno]{amsart}
\usepackage{amssymb}
\usepackage{txfonts}
\usepackage{bbm}
\usepackage{cases}
\usepackage{amsmath}
\usepackage{graphicx}
\usepackage{mathrsfs}
\usepackage{stmaryrd}
\usepackage{amsfonts}
\usepackage{enumerate,amsmath,amssymb,amsthm}

\numberwithin{equation}{section}

\newcommand{\be}{\begin{eqnarray}}
\newcommand{\ee}{\end{eqnarray}}
\newcommand{\ce}{\begin{eqnarray*}}
\newcommand{\de}{\end{eqnarray*}}
\newtheorem{theorem}{Theorem}[section]
\newtheorem{lemma}[theorem]{Lemma}
\newtheorem{remark}[theorem]{Remark}
\newtheorem{definition}[theorem]{Definition}
\newtheorem{proposition}[theorem]{Proposition}
\newtheorem{Examples}[theorem]{Example}
\newtheorem{corollary}[theorem]{Corollary}

\def\v{{\mathbf{v}}}
\def\e{{\mathrm{e}}}
\def\eps{\varepsilon}

\def\p{\partial}

\def\[{{\Big[}}
\def\]{{\Big]}}
\def\<{{\langle}}
\def\>{{\rangle}}
\def\({{\Big(}}
\def\){{\Big)}}

\def\bx{{\mathbf{x}}}

\def\sgn{\mbox{\rm sgn}}

\def\dif{{\mathord{{\rm d}}}}
\def\dis{{\mathord{{\rm \bf d}}}}

\def\no{\nonumber}
\def\={&\!\!=\!\!&}
\def\bt{\begin{theorem}}
\def\et{\end{theorem}}
\def\bl{\begin{lemma}}
\def\el{\end{lemma}}
\def\br{\begin{remark}}
\def\er{\end{remark}}

\def\bd{\begin{definition}}
\def\ed{\end{definition}}
\def\bp{\begin{proposition}}
\def\ep{\end{proposition}}
\def\bc{\begin{corollary}}
\def\ec{\end{corollary}}
\def\bx{\begin{Examples}}
\def\ex{\end{Examples}}

\def\cF{{\mathcal F}}

\def\cP{{\mathcal P}}

\def\mD{{\mathbb D}}
\def\mE{{\mathbb E}}

\def\mH{{\mathbb H}}
\def\mI{{\mathbb I}}

\def\mL{{\mathbb L}}
\def\mM{{\mathbb M}}
\def\mN{{\mathbb N}}

\def\mP{{\mathbb P}}

\def\mR{{\mathbb R}}
\def\mS{{\mathbb S}}

\def\mV{{\mathbb V}}
\def\mW{{\mathbb W}}

\def\sA{{\mathscr A}}
\def\sB{{\mathscr B}}

\def\sE{{\mathscr E}}
\def\sF{{\mathscr F}}
\def\sG{{\mathscr G}}

\def\sI{{\mathscr I}}

\def\geq{\geqslant}
\def\leq{\leqslant}

\def\div{\mathord{{\rm div}}}

\allowdisplaybreaks

\begin{document}

\bf\title{Regularity of density for SDEs driven by degenerate\\
 L\'evy noises}

\date{}
\author{Yulin Song and Xicheng Zhang}

\address{Yulin Song: School of Mathematical Sciences, Beijing Normal University, Beijing, {\rm 100875}, P.R.China\\ Email: songyl@amss.ac.cn}
 
\address{Xicheng Zhang: School of Mathematics and Statistics, Wuhan University, Wuhan, Hubei 430072, P.R.China\\ Email: XichengZhang@gmail.com}

\dedicatory{}

\begin{abstract}
By using Bismut's approach about the Malliavin calculus with jumps, we study the regularity of the distributional density for SDEs driven by degenerate 
additive L\'evy noises. Under full H\"ormander's conditions, we prove the existence of distributional density and the weak continuity in the first variable of the distributional density.
Under the uniform first order Lie's bracket condition, we also prove the smoothness of the density.
\end{abstract}

\maketitle

\rm

\section{Introduction}

Consider the following stochastic differential equation (abbreviated as SDE) in $\mR^d$:
\begin{align}
\dif X_t=b(X_t)\dif t+A_1\dif W_t+ A_2\dif L_t,\ \ X_0=x,\label{SDE}
\end{align}
where $b: \mR^d\to \mR^d$ is a smooth vector field, $A_1$ and $A_2$ are two constant $d\times d$-matrices,
$W_t$ is a $d$-dimensional standard Brownian motion and $L_t$ is a purely jump $d$-dimensional L\'evy process 
with L\'evy measure $\nu(\dif z)$. Let $\Gamma_0:=\{z\in\mR^d: 0<|z|<1\}$. Throughout this work, 
we assume that $\frac{\nu(\dif z)}{\dif z}|_{\Gamma_0}=\kappa(z)$ satisfies the following conditions: for some $\alpha\in(0,2)$ and $m\in\mN$, 
\begin{enumerate}[{\bf (H$^\alpha_m$)}]
\item $\kappa\in C^m(\Gamma_0;(0,\infty))$ is symmetric (i.e. $\kappa(-z)=\kappa(z)$) and satisfies the following Orey's order condition
(cf. \cite[Proposition 28.3]{Sa}): 
\begin{align}
\lim_{\eps\downarrow 0}\eps^{\alpha-2}\int_{|z|\leq \eps}|z|^2\kappa(z)\dif z=:c_1>0,\label{ET11}
\end{align}
and bounded condition: for $j=1,\cdots, m$ and some $C_j>0$,
\begin{align}
|\nabla^j\log\kappa(z)|\leq C_j|z|^{-j}, \ \ z\in \Gamma_0.\label{ET01}
\end{align}
\end{enumerate}
For example, if $\kappa(z)=a(z)|z|^{-d-\alpha}$ with 
$$
a(z)=a(-z), \ \ 0<a_0\leq a(z)\leq a_1,\ \ |\nabla^j a(z)|\leq a_2,\ j=1,\cdots,m,
$$ 
then {\bf (H$^\alpha_m$)} hold. Notice that the generator of SDE (\ref{SDE}) is given by
$$
\sA f(x):=\frac{1}{2}\nabla^2_{A_1A_1^*}f(x)+\mbox{p.v.}\int_{\mR^d}(f(x+A_2z)-f(x))\nu(\dif z)+b(x)\cdot\nabla f(x),
$$
where $A^*_1$ stands for the transpose of $A_1$, and
$$
\nabla^2_{A_1A_1^*}f:=\sum_{i,j=1}^d(A_1A_1^*)_{ij}\p^2_{ij}f,
$$
and p.v. stands for the Cauchy's principle value.

It is well-known that when $b$ is Lipschitz continuous, SDE (\ref{SDE}) has a unique solution $X_t=X_t(x)$.
The aim of this work is to investigate the regularity of the distributional density of $X_t(x)$ under H\"ormander's conditions.
Let $B_0:=\mI_{d\times d}$ be the identity matrix and define for $n\in\mN$,
$$
B_n(x):=b(x)\cdot\nabla B_{n-1}(x)-\nabla b(x)\cdot B_{n-1}(x)+\tfrac{1}{2}\nabla^2_{A_1A_1^*}B_{n-1}(x).
$$ 
Here and below, $(\nabla b(x))_{ij}:=\p_j b^i(x)$. 
Our first main result is about the existence and weak continuity of the 
distribution density for SDE (\ref{SDE}) under full H\"ormander's condition.
\bt\label{Main1}
Assume that {\bf (H$^\alpha_1$)} holds and $b$ is smooth and has bounded derivatives of all orders, and for any $x\in\mR^d$ and some $n=n(x)\in\mN$,
\begin{align}
\mathrm{Rank}[A_1, B_1(x)A_1,\cdots, B_n(x)A_1, A_2, B_1(x)A_2,\cdots, B_n(x)A_2]=d.\label{PT4}
\end{align}
Then $X_t(x)$ admits a density $\rho_t(x,y)$ with respect to the Lebesgue measure so that for any bounded measurable function $f$,
\begin{align}
x\mapsto \cP_t f(x):=\mE f(X_t(x))=\int_{\mR^d} f(y)\rho_t(x,y)\dif y\mbox{ is continuous.}\label{WQ1}
\end{align}
In particular, the semigroup $(\cP_t)_{t\geq 0}$ has the strong Feller property.
\et
\br
When $A_1=0$ and $b(x)=Bx$ is linear, condition (\ref{PT4}) is called Kalman's rank condition. In this case, the smoothness of 
the density of the corresponding Ornstein-Uhlenbeck process has been studied in \cite{Pr-Za, Ku00}.
\er
About the smoothness of the density, we have the following partial result.
\bt\label{Main2}
Assume that {\bf (H$^\alpha_m$)} holds for some $m\in\mN$, and
$b$ is smooth and has bounded derivatives of all orders, and
\begin{align}
\inf_{x\in\mR^d}\inf_{|u|=1}\Big(|uA_1|^2+|uB_1(x)A_1|^2+|uA_2|^2+|uB_1(x)A_2|^2\Big)=:c_2>0.\label{Con}
\end{align}
Then for any $k,n\in\{0\}\cup\mN$ with $1\leq k+n\leq m$, there are $\gamma_{k,n}>0$
and $C=C(k,n)>0$ such that for all $f\in C^\infty_b(\mR^d)$ and $t\in(0,1)$,
\begin{align}
\sup_{x\in\mR^d}|\nabla^k\mE( (\nabla^nf)(X_t(x)))|\leq C\|f\|_\infty t^{-\gamma_{k,n}},
\end{align}
where $\nabla^k$ denotes the $k$-order gradient operator. In particular, if $m=\infty$, then $X_t(x)$ admits a smooth density $\rho_t(x,y)$ so that
\begin{align}
(x,y)\mapsto\rho_t(x,y)\in C^\infty_b(\mR^d\times\mR^d),\ \ t>0.\label{PT5}
\end{align}
\et

In the continuous diffusion case (i.e. $A_2=0$ and $A_1=A_1(x)$), under H\"ormander's conditions, Malliavin \cite{Ma} proved that SDE (\ref{SDE})
has a smooth density by using the stochastic calculus of variations (nowadays, it is also called the Malliavin calculus, and a systematic introduction about the Malliavin calculus
is refereed to the book \cite{Nu}). Since the pioneering work of \cite{Ma}, 
there are many works devoting to extend the Malliavin's theory to the jump case (cf. \cite{Bi0, Bi-Gr-Ja, Pi, Is-Ku} etc.).
However, unlike the case of continuous Brownian functionals, there does not exist a unified treatment for Poisson functionals since the canonical Poisson space has a nonlinear structure.
We mention that Bismut's approach is based on the Girsanov's transformation (cf. \cite{Bi0}), while Picard's approach is to use the difference operator to establish an integration
by parts formula (cf. \cite{Pi}).

When $A_1=0$ and $\kappa(z)= c |z|^{-d-\alpha}$, Theorems \ref{Main1} and \ref{Main2} have been proved in \cite{Zh2} and \cite{Do-Pe-So-Zh} 
by using the Malliavin calculus for subordinated Brownian motions (cf. \cite{Ku}). About the smoothness of distributional density for degenerate
SDEs driven by purely jump noises, Takeuchi \cite{Ta}, Cass \cite{Ca} and Kunita \cite{Ku0} have already studied this problem under different H\"ormander's conditions.
However, their results do not cover the present general case (see also \cite{Zh3,Zh4,Xu} for some related works).
Compared with \cite{Zh2} and \cite{Do-Pe-So-Zh}, in this work we shall use Bismut's approach to prove Theorems \ref{Main1} and \ref{Main2}, and
need to assume that the L\'evy measure is absolutely continuous with respect to the Lebesgue measure.  
It is noticed that in \cite{Do-Pe-So-Zh}, the L\'evy measure can be singular and the drift is allowed to be
arbitrarily growth, which cannot be dealt with in the current settings.

In the proof of our main theorems, one of the difficulties we are facing is the infinity of the moments of $L_t$. To overcome this difficulty, 
we consider two independent L\'evy processes $L^0_t$ and $L^1_t$ with L\'evy measures $\nu_0(\dif z):=1_{|z|<1}\kappa(z)\dif z$ and $\nu_1(\dif z):=1_{|z|\geq 1}\nu(\dif z)$ respectively.
Clearly,
$$
L_t \mbox{ has the same law as $L^0_t+L^1_t$}.
$$
Notice that $L^1_t$ is a compound Poisson process. Let $0=:\tau_0<\tau_1<\tau_2<\cdots<\tau_n<\cdots$ be the jump time of $L^1_t$.
It is well-known that
$$
\sE:=\{\tau_n-\tau_{n-1}, n\in\mN\},\ \ \sG:=\{\Delta L^1_{\tau_n}:=L_{\tau_n}-L_{\tau_n-},n\in\mN\}
$$
are two independent families of i.i.d.
Let $\hbar$ be a c\'adl\'ag purely discontinuous $\mR^d$-valued
function with finite many jumps and $\hbar_0=0$. Following the argument of \cite[Subsection 3.3]{Zh2}, we consider the following SDE:
$$
\tilde X_t(x;\hbar)=x+\int^t_0b(\tilde X_s(x;\hbar))\dif s+A_1W_t+A_2L^0_t+A_2\hbar_t.
$$
Clearly,
$$
X_t(x)\stackrel{(d)}{=}\tilde X_t(x;L^1_\cdot).
$$
If we write
$$
\cP_t f(x):=\mE f(X_t(x)),\ \ \tilde\cP_t f(x):=\mE f(\tilde X_t(x;0)),
$$
then we have (see \cite[(3.19)]{Zh2})
\begin{align}
\cP_t f(x)=\sum_{n=0}^\infty\mE\left(\tilde\cP_{\tau_1}\cdots\vartheta_{A_2\Delta L^1_{\tau_{n-1}}}
\tilde\cP_{\tau_n-\tau_{n-1}}\vartheta_{A_2\Delta L^1_{\tau_n}}\tilde\cP_{t-\tau_n}f(x); \tau_n<t\leq\tau_{n+1}\right),\label{For1}
\end{align}
where for a function $g(x)$ and $y\in\mR^d$,
$$
\vartheta_y g(x):=g(x+y).
$$
Basing on (\ref{For1}) and as in \cite[Subsection 3.3]{Zh2}, it suffices to prove Theorems \ref{Main1} and \ref{Main2} for $\tilde X_t(x;0)$, that is, we only need to consider the SDE 
(\ref{SDE}) driven by $W_t$ and $L^0_t$.

This paper is organised as follows: in Section 2, we recall the Bismut's approach about the Malliavin calculus with jumps.  In \cite{Bi-Gr-Ja},
Bichteler, Gravereaux and Jacod have already systematically introduced it, however, the $\alpha$-stable like noise does not fall into their framework.
Thus, we have to extend the integration by parts formula to the more general class of L\'evy measures. Moreover, we also prove a Kusuoka-Stroock's formula for Poisson stochastic integrals.
In Section 3, we introduce the reduced Malliavin matrix for SDE (\ref{SDE}) used in the Bismut's approach (cf. \cite{Bi-Gr-Ja}), 
and also give some necessary estimates.
In Sections 4 and 5, we shall prove Theorems \ref{Main1} and \ref{Main2}.

Convention: The letter $C$ or $c$ with or without subscripts will denote an unimportant constant, whose value may be different in different places.

\section{Revisit of Bismut's approach to the Malliavin calculus with jumps}

Let $\Gamma\subset\mR^d$ be an open set containing the original point. Let us define
\begin{align}
\Gamma_0:=\Gamma\setminus\{0\},\ \ \varrho(z):=1\vee\dis(z,\Gamma^c_0)^{-1},\label{Rho}
\end{align}
where $\dis(z,\Gamma^c_0)$ is the distance of $z$ to the complement of $\Gamma_0$. 
Let $\Omega$ be the canonical space of all points $\omega=(w,\mu)$, where 
\begin{itemize}
\item $w: [0,1]\to\mR^d$ is a continuous function with $w(0)=0$;
\item $\mu$ is an integer-valued measure on $[0,1]\times\Gamma_0$ with $\mu(A)<+\infty$ for any compact set $A\subset[0,1]\times\Gamma_0$.
\end{itemize} 
Define the canonical process on $\Omega$ as follows: for $\omega=(w,\mu)$,
$$
W_t(\omega):=w(t),\ \ \ N(\omega; \dif t,\dif z):=\mu(\omega; \dif t,\dif z):=\mu(\dif t,\dif z).
$$
Let $(\sF_t)_{t\in[0,1]}$ be the smallest right-continuous filtration on $\Omega$ such that $W$ and $N$ are optional. 
In the following, we write $\sF:=\sF_1$, and endow $(\Omega,\sF)$ with the unique probability measure $\mP$ such that 
\begin{itemize}
\item $W$ is a standard $d$-dimensional Brownian motion;
\item $N$ is a Poisson random measure with intensity $\dif t\nu(\dif z)$, where $\nu(\dif z)=\kappa(z)\dif z$ with
\begin{align}
\kappa\in C^1(\Gamma_0;(0,\infty)),\  \int_{\Gamma_0}(1\wedge|z|^2)\kappa(z)\dif z<+\infty,\ \ |\nabla\log\kappa(z)|\leq C\varrho(z),\label{ET1}
\end{align}
where $\varrho(z)$ is defined by (\ref{Rho});
\item $W$ and $N$ are independent.
\end{itemize}
In the following we write
$$
\tilde N(\dif t,\dif z):=N(\dif t,\dif z)-\dif t\nu(\dif z).
$$
\subsection{Functions spaces}
Let $p\geq 1$. We introduce the following spaces for later use.
\begin{itemize}
\item $\mL^1_p$: The space of all predictable processes: $\xi:\Omega\times[0,1]\times\Gamma_0\to\mR^k$ with finite norm:
$$
\|\xi\|_{\mL^1_p}:=\left[\mE\left(\int^1_0\!\!\!\int_{\Gamma_0}|\xi(s,z)|\nu(\dif z)\dif s\right)^p\right]^{\frac{1}{p}}
+\left[\mE\int^1_0\!\!\!\int_{\Gamma_0}|\xi(s,z)|^p\nu(\dif z)\dif s\right]^{\frac{1}{p}}<\infty.
$$
\item $\mL^2_p$: The space of all predictable processes: $\xi:\Omega\times[0,1]\times\Gamma_0\to\mR^k$ with finite norm:
$$
\|\xi\|_{\mL^2_p}:=\left[\mE\left(\int^1_0\!\!\!\int_{\Gamma_0}|\xi(s,z)|^2\nu(\dif z)\dif s\right)^{\frac{p}{2}}\right]^{\frac{1}{p}}
+\left[\mE\int^1_0\!\!\!\int_{\Gamma_0}|\xi(s,z)|^p\nu(\dif z)\dif s\right]^{\frac{1}{p}}<\infty.
$$
\item $\mH_p$: The space of all measurable adapted processes $h:\Omega\times[0,1]\to\mR^d$ with finite norm:
$$
\|h\|_{\mH_p}:=\left[\mE\left(\int^1_0|h(s)|^2\dif s\right)^{\frac{p}{2}}\right]^{\frac{1}{p}}<+\infty.
$$
\item $\mV_p$:  The space of all predictable processes $\v: \Omega\times[0,1]\times\Gamma_0\to\mR^d$ with finite norm:
$$
\|\v\|_{\mV_p}:=\|\nabla\v\|_{\mL^1_p}+\|\v\varrho\|_{\mL^1_p}<\infty,
$$
where $\varrho(z)$ is defined by (\ref{Rho}). Below we shall write
$$
\mH_{\infty-}:=\cap_{p\geq 1}\mH_p,\ \ \mV_{\infty-}:=\cap_{p\geq 1}\mV_p.
$$
\item $\mH_0$:  The space of all bounded measurable adapted  processes $h:\Omega\times[0,1]\to\mR^d$.
\item $\mV_0$:  The space of all predictable processes $\v: \Omega\times[0,1]\times\Gamma_0\to\mR^d$ 
with the following properties: (i) $\v$ and $\nabla_z \v$ are bounded;
(ii) there exists a compact subset $U\subset \Gamma_0$ such that
$$
\v(t,z)=0,\ \ \forall z\notin U.
$$
\end{itemize}
\br\label{Re21}
For $\xi\in\mL^1_p$, if there is a compact subset $U\subset\Gamma_0$ such that $\xi(s,z)=0$ for all $z\notin U$, then in view of $\kappa\in C^1(\Gamma_0;(0,\infty))$,
$$
\|\xi\|_{\mL^1_p}^p\asymp \mE\left(\int^1_0\!\!\!\int_{U}|\xi(s,z)|^p\dif z\dif s\right)
=\mE\left(\int^1_0\!\!\!\int_{\mR^d}|\xi(s,z)|^p\dif z\dif s\right),
$$
where $\asymp$ means that both sides are comparable up to a constant (depending only on $U,\kappa,p,d$).
\er
\bl\label{Le1}
\begin{enumerate}[(i)]
\item For any $p\geq 1$, the spaces $(\mH_p,\|\cdot\|_{\mH_p})$ and $(\mV_p,\|\cdot\|_{\mV_p})$ are Banach spaces.
\item For any $p_2>p_1\geq 1$, $\mH_{p_2}\subset\mH_{p_1}$ and $\mV_{p_2}\subset\mV_{p_1}$.
\item For any $p\geq 1$, $\mV_0$ (resp. $\mH_0$) is dense in $\mV_p$ (resp.  $\mH_p$).
\end{enumerate}
\el
\begin{proof}
(i) and (ii) are obvious.
\\
\\
(iii) We only prove the density of $\mV_0$ in $\mV_p$, i.e., for each $\v\in\mV_p$, there exists a sequence $\v_n\in\mV_0$ such that 
$$
\lim_{n\to\infty}\|\v_n-\v\|_{\mV_p}=0.
$$
We shall construct the approximation by three steps.\\
\\
(1) For $\eps\in(0,1)$, define
$$
\Gamma_\eps:=\Big\{z\in\mR^d: \dis(z,\Gamma^c_0)>\eps\Big\}.
$$
Let $\chi_\eps:\mR^d\to[0,1]$ be a smooth function with
\begin{align}
\chi_\eps|_{\Gamma_{2\eps}}=1,\ \ \chi_\eps|_{\Gamma^c_\eps}=0,\ \ \|\nabla\chi_\eps\|_\infty\leq C/\eps.\label{CR2}
\end{align}
For $R>1$, let $\chi_R:\mR^d\to[0,1]$ be a smooth function with
\begin{align}
\chi_R(z)=1,\ \ |z|\leq R,\ \  \chi_R(z)=0,\ |z|>2R,\ \ \|\nabla\chi_R\|_\infty\leq C/R.\label{CR1}
\end{align}
Let us define
\begin{align}
\v_{\eps,R}(s,z)=\v(s,z)\chi_\eps(z)\chi_R(z).\label{PT2}
\end{align}
Notice that for $\eps\in(0,1)$ and $R>1$,
\begin{align}
|\nabla\v_{\eps,R}(s,z)-\nabla\v(s,z)|&\leq C\Big(\eps^{-1}1_{z\in\Gamma_\eps\setminus\Gamma_{2\eps}}+R^{-1}1_{R<|z|<2R}\Big)|\v(s,z)|
+\Big(1_{z\in\Gamma^c_{2\eps}}+1_{|z|>R}\Big)|\nabla\v(s,z)|\no\\
&\leq C\varrho(z)\Big(1_{z\in\Gamma^c_{2\eps}}+1_{|z|>R}\Big)|\v(s,z)|+\Big(1_{z\in\Gamma^c_{2\eps}}+1_{|z|>R}\Big)|\nabla\v(s,z)|,\label{BT5}
\end{align}
where $\varrho(z)$ is defined by (\ref{Rho}). By the dominated convergence theorem, we have
$$
\lim_{\eps\downarrow 0, R\uparrow\infty}\|\v_{\eps,R}-\v\|_{\mV_p}=0.
$$
(2) Next we can assume that for some compact set $U\subset\Gamma_0$,
\begin{align}
\v(s,z)=0,\ \ z\notin U.\label{ER5}
\end{align}
Let $\varphi:\mR^d\to[0,1]$ be a smooth function with
$$
\varphi(x)=1,\ \ |x|\leq 1,\ \ \varphi(x)=0,\ \ |x|\geq 2,\ \ \int_{\mR^d}\varphi(x)\dif x=1.
$$
For $\delta\in(0,1)$, set $\varphi_\delta(x):=\delta^{-d}\varphi(\delta^{-1}x)$ and
\begin{align}
\v_\delta(s,z):=\int_{\mR^d}\v(s,z)\varphi_\delta(x-z)\dif z.\label{Phi}
\end{align}
By (\ref{ER5}) and Remark \ref{Re21}, it is easy to see that
\begin{align}
\|\v\|^p_{\mV_p}\asymp\mE\left(\int^1_0\!\!\!\int_{U}(|\v|+|\nabla_z\v|)^p(s,z)\dif z\dif s\right)
=\mE\left(\int^1_0\!\!\!\int_{\mR^d}(|\v|+|\nabla_z\v|)^p(s,z)\dif z\dif s\right).\label{EW1}
\end{align}
Thus,
$$
\lim_{\delta\downarrow 0}\|\v_\delta-\v\|_{\mV_p}=0.
$$
(3) Lastly we assume that $\v$ is smooth in $z$ and satisfies (\ref{ER5}). 
For $R>1$, we construct $\v_R(s,z)$ as follows:
$$
\v_R(\omega,s,z):=\v(\omega,s,z)\cdot 1_{\|\v(\omega,s,\cdot)\|_\infty\leq R}\cdot 1_{\|\nabla\v(\omega,s,\cdot)\|_\infty\leq R}.
$$
Clearly,
$$
\v_R\in\mV_0.
$$
By (\ref{EW1}) and the dominated convergence theorem, we have
$$
\lim_{R\to\infty}\|\v_{R}-\v\|_{\mV_p}=0.
$$
The proof is complete.
\end{proof}
\subsection{Girsanov's theorem}
We need the following Burkholder's inequality.
\bl\label{Le2}
\begin{enumerate}[(i)]
\item For any $p>1$, there is a constant $C_p>0$ such that for any $\xi\in \mL^1_p$,
\begin{align}
\mE\left(\sup_{t\in[0,1]}\left|\int^t_0\!\!\!\int_{\Gamma_0}\xi(s,z)N(\dif s,\dif z)\right|^p\right)\leq C_p\|\xi\|^p_{\mL^1_p}.\label{BT1}
\end{align}
\item For any $p\geq 2$, there is a constant $C_p>0$ such that for any $\xi\in \mL^2_p$,
\begin{align}
\mE\left(\sup_{t\in[0,1]}\left|\int^t_0\!\!\!\int_{\Gamma_0}\xi(s,z)\tilde N(\dif s,\dif z)\right|^p\right)\leq C_p\|\xi\|^p_{\mL^2_p}.\label{BT2}
\end{align}
\end{enumerate}
\el
\begin{proof}
(i) Let us write 
\begin{align}
M_t:=\int^t_0\!\!\!\int_{\Gamma_0}\xi(s,z) \tilde N(\dif s,\dif z)
=\int^t_0\!\!\!\int_{\Gamma_0}\xi(s,z)N(\dif s,\dif z)-\int^t_0\!\!\!\int_{\Gamma_0}\xi(s,z) \nu(\dif z)\dif s.\label{BT3}
\end{align}
For $p>1$, by It\^o's formula, we have
\begin{align*}
\mE|M_t|^p=\mE\left(\int^t_0\!\!\!\int_{\Gamma_0}(|M_{s-}+\xi(s,z)|^p-|M_{s-}|^p-p\xi(s,z)\sgn(M_{s-}) |M_{s-}|^{p-1}) \nu(\dif z)\dif s\right).
\end{align*}
By Doob's maximal inequality and Young's inequality, we further have
\begin{align*}
\mE\left(\sup_{t\in[0,1]}|M_t|^p\right)&\leq C_p\mE |M_1|^p\leq
C_p\mE\left(\int^1_0\!\!\!\int_{\Gamma_0}|\xi(s,z)|(|M_{s-}|+|\xi(s,z)|)^{p-1}\nu(\dif z)\dif s\right)\\
&\leq C_p\mE\left(\sup_{s\in[0,1]}|M_s|^{p-1}\int^1_0\!\!\!\int_{\Gamma_0}|\xi(s,z)|\nu(\dif z)\dif s\right)+C_p\|\xi\|^p_{\mL^1_p}\\
&\leq \frac{1}{2}\mE\left(\sup_{s\in[0,1]}|M_s|^p\right)+C_p\|\xi\|^p_{\mL^1_p},
\end{align*}
which together with (\ref{BT3}) gives (\ref{BT1}).
\\
\\
(ii) As above, for $p\geq 2$, by Taylor's expansion, we have
\begin{align*}
\mE\left(\sup_{t\in[0,1]}|M_t|^p\right)&\leq C_p\mE\left(\int^1_0\!\!\!\int_{\Gamma_0}|\xi(s,z)|^2(|M_{s-}|+|\xi(s,z)|)^{p-2}\nu(\dif z)\dif s\right)\\
&\leq C_p\mE\left(\sup_{s\in[0,1]}|M_s|^{p-2}\int^1_0\!\!\!\int_{\Gamma_0}|\xi(s,z)|^2\nu(\dif z)\dif s\right)+C_p\|\xi\|^p_{\mL^2_p}\\
&\leq \frac{1}{2}\mE\left(\sup_{s\in[0,1]}|M_s|^p\right)+C_p\|\xi\|^p_{\mL^2_p},
\end{align*}
which in turn gives (\ref{BT2}).
\end{proof}
For $\v\in\mV_0$ and $\eps>0$, define
$$
\gamma_\eps(t,z):=\det\left(I+\eps\nabla_z\v(t,z)\right)\frac{\kappa(z+\eps\v(t,z))}{\kappa(z)}.
$$
The following lemma is easy. 
\bl\label{Le19}
For any $\v\in\mV_0$ with compact support $U\subset\Gamma_0$ with respect to $z$, 
there exist an $\eps_0>0$ and a constant $C>0$ such that for any $\eps\in(0,\eps_0)$ and all $t,z$,
\begin{align}
|\gamma_\eps(t,z)-1|\leq C\eps 1_{U}(z).\label{ER2}
\end{align}
Moreover, we have
\begin{align}
\frac{\dif\gamma_\eps(t,z)}{\dif\eps}\Big|_{\eps=0}=\div \v(t,z)+\<\nabla\log\kappa(z),\v(t,z)\>_{\mR^d}=\frac{\div(\kappa\v)(t,z)}{\kappa(z)}.\label{ER3}
\end{align}
\el
\begin{proof}
Since $\v(t,z)=0$ for $z\notin U$, we have
$$
\gamma_\eps(t,z)=1,\ \ \forall z\notin U.
$$
For any $z\in U$, since $\v$ and $\nabla_z\v$ are bounded, we have
\begin{align*}
|\gamma_\eps(t,z)-1|&\leq |\det\left(I+\eps\nabla_z\v(t,z)\right)|\left|\frac{\kappa(z+\eps\v(t,z))}{\kappa(z)}-1\right|
+|\det\left(I+\eps\nabla_z\v(t,z)\right)-1|\\
&\leq \frac{C}{\inf_{z\in U}\kappa(z)}|\kappa(z+\eps\v(t,z))-\kappa(z)|+C\eps,
\end{align*}
which gives the desired estimate (\ref{ER2}) by the compactness of $U$ and $\kappa\in C^1(\Gamma_0; (0,\infty))$. As for (\ref{ER3}), it follows by a direct calculation.
\end{proof}

For $p\geq 1$ and  $\Theta:=(h,\v)\in\mH_p\times\mV_p$,  we write
\begin{align}
\div\Theta:=-\int^1_0\<h(s),\dif W_s\>_{\mR^d}+\int^1_0\!\!\!\int_{\Gamma_0}\frac{\div(\kappa\v)(s,z)}{\kappa(z)}\tilde N(\dif s,\dif z).\label{ER1}
\end{align}
By Burkholder's inequality and (\ref{ET1}), we have
\begin{align}
\mE|\div\Theta|^p\leq C\Big(\|h\|_{\mH_p}^p+\|\v\|_{\mV_p}^p\Big).\label{214}
\end{align}
Let $Q^\eps_t$ solve the following SDE:
\begin{align}
Q^\eps_t=1-\eps\int^t_0Q^\eps_s\<h_s,\dif W_s\>_{\mR^d}+\int^t_0\!\!\!\int_{\Gamma_0}Q^\eps_{s-}(\gamma_\eps(s,z)-1)\tilde N(\dif s,\dif z),\label{EQ1}
\end{align}
whose solution is explicitly given by the Doleans-Dade's formula:
\begin{align*}
Q^\eps_t&=\exp\Bigg\{\int^t_0\!\!\!\int_{\Gamma_0}\log \gamma_\eps(s,z)N(\dif s,\dif z)-\int^t_0\!\!\!\int_{\Gamma_0}(\gamma_\eps(s,z)-1)\nu(\dif z)\dif s\\
&\qquad\qquad-\eps\int^t_0\<h_s,\dif W_s\>_{\mR^d}-\frac{\eps^2}{2}\int^t_0|h_s|^2\dif s\Bigg\}.
\end{align*}
\bl
If $\Theta=(h,\v)\in\mH_0\times\mV_0$, then $Q^\eps_t$ is a nonnegative martingale and for any $p\geq 2$,
\begin{align}
\lim_{\eps\downarrow 0}\mE\left|\frac{Q^\eps_1-1}{\eps}-\div\Theta\right|^p=0.\label{ER4}
\end{align}
\el
\begin{proof}
For any $p\geq 2$, by (\ref{EQ1}), (\ref{ER2})  and (\ref{BT1}), we have
\begin{align*}
\mE|Q^\eps_t|^p&\leq C+C\eps^p\int^t_0\mE (|Q^\eps_s|^p|h_s|^p)\dif s+
\int^t_0\!\!\!\int_{\Gamma_0}\mE(|Q^\eps_{s-}(\gamma_\eps(s,z)-1)|^p)\nu(\dif z)\dif s\\
&\leq C+C\eps^p\left(\|h\|_\infty^p+\nu(U)\right)\int^t_0\mE|Q^\eps_s|^p\dif s,
\end{align*}
which gives
\begin{align}
\sup_{\eps\in(0,1)}\sup_{t\in[0,1]}\mE|Q^\eps_t|^p<+\infty.\label{BT4}
\end{align}
From this and (\ref{EQ1}), one sees that $Q^\eps_t$ is a nonnegative martingale and $\mE Q^\eps_t=1$.

For (\ref{ER4}), by equation (\ref{EQ1}) and (\ref{BT4}),  we have
$$
\lim_{\eps\downarrow 0}\sup_{t\in[0,1]}\mE|Q^\eps_t-1|^p=0,
$$
and
\begin{align*}
\frac{Q^\eps_t-1}{\eps}-\div\Theta
&=\int^1_0(Q^\eps_s-1)\<h_s,\dif W_s\>_{\mR^d}
+\frac{1}{\eps}\int^1_0\!\!\!\int_{\Gamma_0}(Q^\eps_{s-}-1)(\gamma_\eps(s,z)-1)\tilde N(\dif s,\dif z)\\
&\quad+\int^1_0\!\!\!\int_{\Gamma_0}\left(\frac{\gamma_\eps(s,z)-1}{\eps}-\frac{\div(\kappa\v)(s,z)}{\kappa(z)}\right)\tilde N(\dif s,\dif z).
\end{align*}
Thus, by Burkholder's inequality and Lemma \ref{Le19}, we obtain (\ref{ER4}).
\end{proof}
For $\Theta=(h,\v)\in\mH_0\times\mV_0$ and $\eps>0$, define
$$
W^\eps_t:=W_t+\eps\int^t_0h(s)\dif s,\ \ N^\eps((0,t]\times E):=\int^t_0\!\!\!\int_{\Gamma_0}1_E(z+\eps\v(s,z))N(\dif s,\dif z).
$$
Then the map
\begin{align}
\Theta^\eps: (W, N)\mapsto (W^\eps, N^\eps)\label{Theta}
\end{align}
defines a transformation from $\Omega$ to $\Omega$.
We have (cf. \cite[p.64, Theorem 6-16]{Bi-Gr-Ja} or \cite[p. 185]{Bi})
\bt\label{Gir}
(Girsanov's theorem) For $\Theta=(h,\v)\in\mH_0\times\mV_0$, there exists an $\eps_0>0$ such that for all $\eps\in(0,\eps_0)$, the law of $(W^\eps, N^\eps)$ under $Q^\eps_1\mP$
is the same as $\mP$, i.e.,
$$
\mP=(Q^\eps_1\mP)\circ(\Theta^\eps)^{-1}.
$$
\et
\subsection{Malliavin derivative operator}
Let $C_{\mathrm{p}}^\infty(\mR^m)$ be the class of all smooth functions on $\mR^m$ which together with all the derivatives have at most polynomial growth. 
Let $\cF C^\infty_{\mathrm{p}}$ be the class of all Wiener-Poisson functionals on $\Omega$ with the following form:
$$
F(\omega)=f(w(h_1),\cdots, w(h_{m_1}), \mu(g_1),\cdots, \mu(g_{m_2})),\ \ \omega=(w,\mu)\in\Omega,
$$
where $f\in C_{\mathrm{p}}^\infty(\mR^{m_1+m_2})$, $h_1,\cdots, h_{m_1}\in\mH_0$ and $g_1,\cdots, g_{m_2}\in\mV_0$ are non-random, and
$$
w(h_i):=\int^1_0\<h_i(s), \dif w(s)\>_{\mR^d},\ \ \mu(g_j):=\int^1_0\!\!\!\int_{\Gamma_0}g_j(s,z)\mu(\dif s,\dif z).
$$
Notice that
$$
\cF C^\infty_{\mathrm{p}}\subset \cap_{p\geq 1}L^p(\Omega,\sF,\mP).
$$
For $\Theta=(h,\v)\in\mH_{\infty-}\times\mV_{\infty-}$, let us define
\begin{align}
D_\Theta F&:=\sum_{i=1}^{m_1}(\p_i f)(\cdot)\int^1_0\<h(s), h_i(s)\>_{\mR^d}\dif s\no\\
&+\sum_{j=1}^{m_2}(\p_{j+m_1} f)(\cdot)\int^1_0\!\!\!\int_{\Gamma_0}\v(s,z)\cdot \nabla_z g_j(s,z)\mu(\dif s,\dif z),\label{PT1}
\end{align}
where ``$(\cdot)$'' stands for $w(h_1),\cdots, w(h_{m_1}), \mu(g_1),\cdots, \mu(g_{m_2})$. 
By H\"older's inequality and (\ref{BT2}), it is easy to see that for any $p\geq 1$,
\begin{align}
D_\Theta F\in L^p\mbox{ and } \ D_\Theta F=\lim_{\eps\to 0}\frac{F\circ\Theta^\eps-F}{\eps}\ \mbox{ in $L^p$},\label{Dif}
\end{align}
where $\Theta^\eps$ is defined by (\ref{Theta}).
Thus, $D_\Theta F$ is well defined, i.e., it does not depend on the representation of $F$.

We have
\bl
Let $\Theta=(h,\v)\in \mH_{\infty-}\times\mV_{\infty-}$ and $\div\Theta$ be defined by (\ref{ER1}). 
\begin{enumerate}[(i)]
\item (Density) $\cF C^\infty_{\mathrm{p}}$ is dense in $L^p:=L^p(\Omega,\sF,\mP)$ for any $p\geq 1$.
\item (Integration by parts) For any $F\in\cF C^\infty_{\mathrm{p}}$, we have
\begin{align}
\mE(D_\Theta F)=-\mE(F \div\Theta).\label{ER8}
\end{align}
\item (Closability) The linear operator $(D_\Theta, \cF C^\infty_{\mathrm{p}})$ is closable in $L^p$ for any $p>1$.
\end{enumerate}
\el
\begin{proof}
(i) is standard by a monotonic argument.
\\
\\
(ii) We first assume $\Theta=(h,\v)\in\mH_0\times\mV_0$. By (\ref{Dif}) and Theorem \ref{Gir}
, we have
\begin{align*}
\mE D_\Theta F&=\lim_{\eps\downarrow 0}\frac{1}{\eps}\mE(F\circ\Theta^\eps-F)
=\lim_{\eps\downarrow 0}\frac{1}{\eps}\mE((1-Q^\eps_1) F\circ\Theta^\eps)=-\mE(F \div\Theta),
\end{align*}
where we have used (\ref{ER4}) in the last step. 
For general $\Theta=(h,\v)\in\mH_{\infty-}\times\mV_{\infty-}$ and $p>2$, 
by Lemma \ref{Le19} there exists a sequence of $\Theta_n=(h_n,\v_n)\in\mH_0\times\mV_0$
such that
$$
\lim_{n\to\infty}(\|h_n-h\|_{\mH_p}+\|\v_n-\v\|_{\mV_p})=0.
$$
By the definition of $D_\Theta F$, it is easy to see that
$$
\lim_{n\to\infty}\mE|D_{\Theta_n}F-D_\Theta F|^2=0.
$$
Moreover, by (\ref{214}) we also have
$$
\lim_{n\to\infty}\mE|\div(\Theta_n -\Theta)|^2\leq \lim_{n\to\infty}(\|h_n-h\|^2_{\mH_2}+\|\v_n-\v\|_{\mV_2}^2)=0.
$$
By taking limits for $\mE(D_{\Theta_n} F)=-\mE(F \div\Theta_n)$, we obtain (\ref{ER8}).
\\
\\
(iii) Fix $p>1$. Let $F_n$ be a sequence in $\cF C^\infty_{\mathrm{p}}$ converging to zero in $L^p$. Suppose that $D_\Theta F_n$ converges to some $\xi$ in $L^p$.
We want to show $\xi=0$. For any $G\in\cF C^\infty_{\mathrm{p}}$, noticing that $F_n G\in\cF C^\infty_{\mathrm{p}}$, 
by H\"older's inequality, we have
\begin{align*}
\mE(G\xi)=\lim_{n\to\infty}\mE(GD_\Theta F_n)&\stackrel{(\ref{Dif})}{=}
\lim_{n\to\infty}\mE(D_\Theta (F_n G))-\lim_{n\to\infty}\mE(F_nD_\Theta G)\\
&\stackrel{(\ref{ER8}) }{=}-\lim_{n\to\infty}\mE(F_n G\div\Theta)=0.
\end{align*}
By (i), we obtain $\xi=0$. The proof is complete.
\end{proof}
\bd
For given $\Theta=(h,\v)\in\mH_{\infty-}\times\mV_{\infty-}$ and $p>1$,  we define the first order Sobolev space $\mW_{\Theta}^{1,p}$ being 
the completion of $\cF C^\infty_{\mathrm{p}}$ in $L^p(\Omega,\sF,\mP)$ with respect to the norm:
$$
\|F\|_{\Theta; 1,p}:=\|F\|_{L^p}+\|D_\Theta F\|_{L^p}.
$$
Clearly,  $\mW_{\Theta}^{1,p_2}\subset\mW_{\Theta}^{1,p_1}$ for $p_2>p_1>1$.
We shall write
$$
\mW_{\Theta}^{1,\infty-}:=\cap_{p>1}\mW_{\Theta}^{1,p}.
$$
\ed
We have the following integration by parts formula.
\bt
Let $\Theta=(h,\v)\in\mH_{\infty-}\times\mV_{\infty-}$ and $p>1$. For any $F\in\mW^{1,p}_\Theta$, we have
\begin{align}
\mE(D_\Theta F)=-\mE(F \div\Theta),\label{ER88}
\end{align}
where $\div\Theta$ is defined by (\ref{ER1}).
\et
\begin{proof}
Let $F_n\in\cF C^\infty_{\mathrm{p}}$ converge to $F$ in $\mW^{1,p}_\Theta$. By (\ref{ER8}) we have
$$
\mE(D_\Theta F_n)=-\mE(F_n \div\Theta).
$$
By taking limits, we obtain (\ref{ER88}).
\end{proof}
Moreover, we also have the following chain rule.
\bp (Chain rule) Let $\Theta=(h,\v)\in\mH_{\infty-}\times\mV_{\infty-}$.
For $m,k\in\mN$, let $F=(F_1,\cdots, F_m)\in (\mW^{1,\infty-}_\Theta)^m$ and $\varphi\in C^\infty_{\mathrm{p}}(\mR^m;\mR^k)$.
Then the composition $\varphi(F)\in (\mW^{1,\infty-}_\Theta)^k$ and 
$$
D_\Theta\varphi (F)=D_\Theta F\cdot\nabla\varphi(F).
$$
\ep
\begin{proof}
Since $\varphi\in C^\infty_{\mathrm{p}}(\mR^m;\mR^k)$, we can assume that for some $r\in\mN$,
\begin{align}
|\nabla\varphi(x)|\leq C(1+|x|^r).\label{KL1}
\end{align}
For any fixed $p>r+1$, let $F_n\in(\cF C^\infty_\mathrm{p})^m$ converge to $F$ in $(\mW^{1,p}_\Theta)^m$.
Since $\varphi (F_n)\in(\cF C^\infty_{\mathrm{p}})^k$,  by (\ref{Dif}) it is easy to see that
$$
D_{\Theta}\varphi (F_n)=D_\Theta F_n\cdot\nabla\varphi(F_n).
$$
For any $q\in(1,\frac{p}{r+1})$, by H\"older's inequality and (\ref{KL1}), we have
\begin{align*}
&\lim_{n\to\infty}\mE|D_\Theta F_n\cdot\nabla\varphi(F_n)-D_\Theta F\cdot\nabla\varphi(F)|^q\\
&\quad\leq C\lim_{n\to\infty}\Big(\mE|D_\Theta F_n|^p\Big)^{\frac{q}{p}}\left(\mE|\nabla\varphi(F_n)-\nabla\varphi(F)|^{\frac{qp}{p-q}}\right)^{\frac{p-q}{p}}\\
&\quad+C\lim_{n\to\infty}\Big(\mE|D_\Theta F_n-D_\Theta F|^p\Big)^{\frac{q}{p}}\left(1+\mE|F|^{\frac{rqp}{p-q}}\right)^{\frac{p-q}{p}}=0,
\end{align*}
and also
$$
\lim_{n\to\infty}\mE|\varphi(F_n)-\varphi(F)|^q\leq C\lim_{n\to\infty}\Big(\mE|F_n-F|^p\Big)^{\frac{q}{p}}
\left(1+\mE|F_n|^{\frac{rqp}{p-q}}+\mE|F|^{\frac{rqp}{p-q}}\right)^{\frac{p-q}{p}}=0.
$$
Thus, by definition we have $\varphi (F)\in(\mW^{1,q}_\Theta)^k$ and
$$
D_{\Theta}\varphi (F)=D_\Theta F\cdot\nabla\varphi(F).
$$
Since $p>r+1$ is arbitrary and $q\in(1,\frac{p}{r+1})$, we obtain $\varphi (F)\in(\mW^{1,\infty-}_\Theta)^k$.
\end{proof}

\subsection{Kusuoka-Stroock's formula}
In this subsection we are about to establish a commutation formula between the gradient and Poisson stochastic integrals. On Wiener space this formula is
given by Kusuoka and Stroock \cite{Ku-St1}. On configuration space similar formula is proven in \cite{Re-Ro-Zh}.
\bp\label{Pr1}
Fix $\Theta=(h,\v)\in\mH_{\infty-}\times\mV_{\infty-}$.  Let $\eta(\omega,s,z):\Omega\times[0,1]\times\Gamma_0\to\mR$ be a measurable map 
and satisfy that for each $(s,z)\in[0,1]\times\Gamma_0$, 
$$\eta(s,z)\in\mW^{1,\infty-}_\Theta,\ \ \eta(s,\cdot)\in C^1(\Gamma_0)
$$ 
and
\begin{align}
\mbox{$s\mapsto \eta(s,z), D_\Theta\eta(s,z),\nabla_z\eta(s,z)$ are left-continuous and $\sF_s$-adapted},
\end{align}
and  for any $p>1$, 
\begin{align}
\mE\left[\sup_{s\in[0,1]}\sup_{z\in\Gamma_0}\left(\frac{|\eta(s,z)|^p+|D_\Theta\eta(s,z)|^p}{(1\wedge|z|)^p}+|\nabla_z\eta(s,z)|^p\right)\right]<+\infty.\label{ET7}
\end{align}
Then $\sI(\eta):=\int^1_0\!\int_{\Gamma_0}\eta(s,z)\tilde N(\dif s,\dif z)\in\mW^{1,\infty-}_\Theta$ and
\begin{align}
D_\Theta \sI(\eta)=\int^1_0\!\!\!\int_{\Gamma_0}D_\Theta\eta(s,z)\tilde N(\dif s,\dif z)+\int^1_0\!\!\!\int_{\Gamma_0}\v(s,z) \cdot\nabla\eta(s,z)N(\dif s,\dif z).\label{For2}
\end{align}
\ep
\begin{proof}
(i) First of all, we assume that $\eta(s,z)=1_{(t_0,t_1]}(s)\eta(z)$, where $\eta(z)$ is $\sF_{t_0}$-measurable, and satisfies (\ref{ET7}) and
\begin{align}
\mbox{$z\mapsto \eta(z)$ has compact support $U\subset\Gamma_0$.}
\end{align}
For $n\in\mN$, let $\mD_n$ be the grid of $\mR^d$ with step $2^{-n}$. For a point $z\in\mR^d$, 
let $\phi_n(z)$ be the left-lower corner point in $\mD_n$ which is closest to $z$. 
For $\eps\in(0,1)$ and $R>1$, let $\chi_\eps$ and $\chi_R$ be defined by (\ref{CR2}) and (\ref{CR1}). 
For $\delta\in(0,1)$, let $\eta_\delta(z)$ be defined as in (\ref{Phi}), and let us define
$$
\eta^{\delta,n}_{\eps,R}(\omega,y):=\chi_\eps(y)\chi_R(y)\int^{y_1}_0\cdots\int^{y_d}_0(\p_{z_1}\cdots\p_{z_d}\eta_\delta)(\omega,\phi_n(z))\dif z_1\cdots\dif z_d.
$$
From this definition, we can write
$$
\eta^{\delta,n}_{\eps,R}(\omega,z)=\sum_{k=1}^m\xi_j(\omega)g_j(z),
$$
where $\xi_j\in\mW^{1,\infty-}_\Theta$ is $\sF_{t_0}$-measurable and $g_j$ is smooth and has support 
$$
U_{\eps,R}:=\Gamma_\eps\cap\{z: |z|\leq 2R\}\subset\Gamma_0.
$$
By definition (\ref{PT1}), it is easy to check that $\sI(\eta^{\delta,n}_{\eps,R}):=\int_U\eta^{\delta,n}_{\eps,R}(z)\tilde N((t_0,t_1],\dif z)\in \mW^{1,\infty-}_\Theta$ and
$$
D_\Theta\sI(\eta^{\delta,n}_{\eps,R})=\int_{U_{\eps,R}}D_\Theta\eta^{\delta,n}_{\eps,R}(z)\tilde N((t_0,t_1],\dif z)
+\int^{t_1}_{t_0}\!\!\!\int_{U_{\eps,R}}\v(s,z)\cdot\nabla\eta^{\delta,n}_{\eps,R}(z)N(\dif s,\dif z).
$$
Thus, for proving (\ref{For2}), by Lemma \ref{Le2} it suffices to prove that for any $p>1$,
\begin{align*}
\lim_{R\to\infty}\lim_{\eps\to 0}\lim_{\delta\to\infty}\lim_{n\to\infty}\Big(\|\eta^{\delta,n}_{\eps,R}-\eta\|_{\mL_p^1}+\|D_\Theta(\eta^{\delta,n}_{\eps,R}-\eta)\|_{\mL_p^1}\Big)=0,
\end{align*}
and\begin{align}
\lim_{R\to\infty}\lim_{\eps\to 0}\lim_{\delta\to\infty}\lim_{n\to\infty}\|\v\cdot\nabla (\eta^{\delta,n}_{\eps,R}-\eta)\|_{\mL_p^1}=0.\label{Lim12}
\end{align}
We only prove the second limit.  The first limit is similar. For fixed $\eps, R$, set $\eta_{\eps,R}:=\chi_\eps\chi_R\eta$. Since for $z\not \in U_{\eps,R}$,
$$
\eta^{\delta,n}_{\eps,R}(z)=\eta_{\eps,R}(z)=0,
$$
by Remark \ref{Re21} and H\"older's inequality, we have
\begin{align}
\lim_{\delta\to 0}\lim_{n\to\infty}\|\v\cdot\nabla (\eta^{\delta,n}_{\eps,R}-\eta_{\eps,R})\|_{\mL_p^1}^p
&\leq C\lim_{\delta\to 0}\lim_{n\to\infty}\mE\int^{t_1}_{t_0}\!\!\!\int_{U_{\eps,R}}|\v(s,z)\cdot\nabla (\eta^{\delta,n}_{\eps,R}-\eta_{\eps,R})(s,z)|^p\dif z\dif s\no\\
&\leq C\lim_{\delta\to 0}\lim_{n\to\infty}\left(\mE\int^{t_1}_{t_0}\!\!\!\int_{U_{\eps,R}}|\nabla (\eta^{\delta,n}_{\eps,R}-\eta_{\eps,R})(s,z)|^{2p}\dif z\dif s\right)^{\frac{1}{2}}=0.\label{Lim13}
\end{align}
On the other hand, since $\eta$ has compact support $U$, by (\ref{ET7}) and the dominated convergence theorem, we have
\begin{align}
\lim_{R\to \infty}\lim_{\eps\to 0}\|\v\cdot\nabla (\eta_{\eps,R}-\eta)\|_{\mL_p^1}^p=
\lim_{R\to \infty}\lim_{\eps\to 0}\|(1-\chi_\eps\chi_R)\v\cdot\nabla \eta\|_{\mL_p^1}^p=0.\label{Lim14}
\end{align}
Combining (\ref{Lim13}) and (\ref{Lim14}), we obtain (\ref{Lim12}).
\\
\\
(ii) Next we assume that for some compact set $U\subset\Gamma_0$,
\begin{align}
\eta(s,z)=0,\ \ z\notin U.\label{Com}
\end{align}
For $n\in\mN$, let $s_k:=k/n$ and define
$$
\eta_n(s,z):=\sum_{k=1}^n1_{(s_{k-1},s_k]}(s)\eta(s_{k-1},z).
$$
In this case, we have
\begin{align*}
\sI(\eta_n )=\sum_{k=1}^n\left(\int_{\Gamma_0}\eta(s_{k-1},z)N((s_{k-1},s_k],\dif z)-\frac{1}{n}\int_{\Gamma_0}\eta(s_{k-1},z)\nu(\dif z)\right)
\end{align*}
By (i), we have
$$
\sI(\eta_n )\in\mW^{1,\infty-}_\Theta
$$
and
\begin{align*}
D_\Theta \sI(\eta_n)
=\int^1_0\!\!\!\int_{\Gamma_0}D_\Theta\eta_n(s,z)\tilde N(\dif s,\dif z)+\int^1_0\!\!\!\int_{\Gamma_0}\v(s,z) \cdot\nabla\eta_n(s,z)N(\dif s,\dif z).
\end{align*}
By Lemma \ref{Le2} and (\ref{Com}), for any $p\geq 2$, we have
\begin{align*}
&\mE\left|\int^1_0\!\!\!\int_{\Gamma_0}(D_\Theta\eta_n(s,z)-D_\Theta\eta(s,z))\tilde N(\dif s,\dif z)\right|^p\\
&\quad\leq C\mE\left(\int^1_0\!\!\!\int_{U}|D_\Theta\eta_n(s,z)-D_\Theta\eta(s,z)|^p\dif z\dif s\right)
\end{align*}
and
\begin{align*}
&\mE\left|\int^1_0\!\!\!\int_{\Gamma_0}\v(s,z) \cdot\nabla(\eta_n-\eta)(s,z)N(\dif s,\dif z)\right|^p\\
&\quad\leq C\mE\left(\int^1_0\!\!\!\int_{U}|\nabla\eta_n(s,z)-\nabla\eta(s,z)|^p~ |\v(s,z)|^p\dif z\dif s\right).
\end{align*}
By the assumptions and the dominated convergence theorem, both of them converges to zero as $n\to\infty$, and we obtain (\ref{For2}).
\\
\\
(iii) We now drop the assumption (\ref{Com}). Define
$$
\eta_{\eps, R}(s,z):=\chi_\eps(z)\chi_R(z)\eta(s,z),
$$
where $\chi_\eps$ and $\chi_R$ are the same as in (\ref{CR2}) and (\ref{CR1}). By (ii), we have
\begin{align*}
D_\Theta \sI(\eta_{\eps, R})
=\int^1_0\!\!\!\int_{\Gamma_0}D_\Theta\eta_{\eps, R}(s,z)\tilde N(\dif s,\dif z)+\int^1_0\!\!\!\int_{\Gamma_0}\v(s,z) \cdot\nabla\eta_{\eps, R}(s,z)N(\dif s,\dif z).
\end{align*}
For proving (\ref{For2}), it suffices to prove that for any $p\geq 2$,
\begin{align*}
&I^{(1)}_{\eps,R}:=\mE\left|\int^1_0\!\!\!\int_{\Gamma_0}(1-\chi_\eps(z)\chi_R(z))D_\Theta\eta(s,z)\tilde N(\dif s,\dif z)\right|^p\to0,\\
&I^{(2)}_{\eps,R}:=\mE\left|\int^1_0\!\!\!\int_{\Gamma_0}\v(s,z) \cdot\nabla(\eta_{\eps, R}-\eta)(s,z)N(\dif s,\dif z)\right|^p\to0,
\end{align*}
as $R\to\infty$ and $\eps\to 0$.
The first limit follows by (\ref{BT2}), (\ref{ET7}) and the dominated convergence theorem. For the second limit, noticing that as in (\ref{BT5}),
$$
|\nabla \eta_{\eps, R}(s,z)-\nabla \eta(s,z)|\leq C\varrho(z)\Big(1_{z\in\Gamma^c_{2\eps}}+1_{|z|>R}\Big)|\eta(s,z)|+\Big(1_{z\in\Gamma^c_{2\eps}}+1_{|z|>R}\Big)|\nabla\eta(s,z)|,
$$
by (\ref{BT1}) we have
\begin{align*}
&\mE\left|\int^1_0\!\!\!\int_{\Gamma_0}\v(s,z) \cdot\nabla(\eta_{\eps, R}-\eta)(s,z)N(\dif s,\dif z)\right|^p\\
&\quad\leq C\mE\left(\sup_{s,z}(|\eta(s,z)|+|\nabla\eta(s,z)|)\int^1_0\!\!\!\int_{\Gamma^c_{2\eps}\cup\{|z|>R\}}|\varrho(z)\v(s,z)|\nu(\dif z)\dif s\right)^p\\
&\quad+C\mE\left(\sup_{s,z}(|\eta(s,z)|+|\nabla\eta(s,z)|)^p\int^1_0\!\!\!\int_{\Gamma_{2\eps}^c\cup\{|z|>R\}}|\varrho(z)\v(s,z)|^p\nu(\dif z)\dif s\right),
\end{align*}
which converges to zero as $\eps\to 0$ and $R\to\infty$. The proof is complete.
\end{proof}

\section{Reduced Malliavin matrix for SDEs driven by L\'evy noises}

As discussed in the introduction, in the remainder of this paper, we shall assume that 
$$
\Gamma_0=\{z\in\mR^d: 0<|z|<1\},
$$ 
and  
$$
\mbox{$\frac{\nu(\dif z)}{\dif z}|_{\Gamma_0}=\kappa(z)$ with $\kappa$ satisfying {\bf (H$^\alpha_1$)}}.
$$
Let $N(\dif t,\dif z)$ be the Poisson random measure associated with $L^0_t$, i.e.,
$$
N((0,t]\times E):=\sum_{s\leq t}1_{E}(L^0_s-L^0_{s-}),\ \ E\in\sB(\Gamma_0).
$$
Since $\nu(\dif z)$ is symmetric, by L\'evy-It\^o's decomposition, we can write
$$
L^0_t=\int^t_0\!\!\!\int_{\Gamma_0}z\tilde N(\dif s,\dif z)=\int^t_0\!\!\!\int_{\Gamma_0}z(N(\dif s,\dif z)-\dif s\nu(\dif z)).
$$
By Proposition \ref{Pr1}, for any $\Theta=(h,\v)\in\mH_{\infty-}\times\mV_{\infty-}$,  
we have $W_t, L^0_t\in\mW^{1,\infty-}_\Theta$ and
\begin{align}
D_\Theta W_t=\int^t_0 h(s)\dif s,\ \ D_\Theta L^0_t= \int^t_0\!\!\!\int_{\Gamma_0}\v(s,z)N(\dif s,\dif z).\label{EG2}
\end{align}
Let $X_t=X_t(x)$ solve the following SDE:
$$
\dif X_t=b(X_t)\dif t+A_1\dif W_t+A_2\dif L^0_t,\ \ X_0=x.
$$
\bp
Assume that $b\in C^1$ has bounded derivative. For fixed $\Theta=(h,\v)\in\mH_{\infty-}\times\mV_{\infty-}$, we have $X_t\in\mW^{1,\infty-}_\Theta$ and
\begin{align}
D_\Theta X_t=\int^t_0\nabla b(X_s) D_\Theta X_s\dif s+A_1 \int^t_0h(s)\dif s+A_2\int^t_0\!\!\!\int_{\Gamma_0}\v(s,z)N(\dif s,\dif z).\label{EW6}
\end{align}
\ep
\begin{proof}
Consider the following Picard's iteration: $X^0_t=x$ and for $n\in\mN$,
$$
X^n_t=x+\int^t_0b(X^{n-1}_s)\dif s+A_1 W_t+ A_2 L^0_t.
$$
It is by now standard to prove that for any $t\geq 0$ and $p\geq 1$,
\begin{align}
\lim_{n\to\infty}\mE|X^n_t-X_t|^p=0.\label{EW10}
\end{align}
Since $\Theta\in\mH_{\infty-}\times\mV_{\infty-}$, by (\ref{EG2}) and the induction, we have $X^n_t\in\mW^{1,\infty-}_\Theta$ and
$$
D_\Theta X^n_t=\int^t_0\nabla b(X^{n-1}_s) D_\Theta X^{n-1}_t\dif s+A_1\int^t_0h(s)\dif s+ A_2\int^t_0\!\!\!\int_{\Gamma_0}\v(s,z)N(\dif s,\dif z).
$$
By Gronwall's inequality, it is easy to prove that for any $T>0$ and $p\geq1$,
$$
\sup_{n\in\mN}\sup_{t\in[0,T]}\mE|D_\Theta X^{n}_t|^p<+\infty.
$$
Let $Y_t$ solve the following SDE:
$$
Y_t=\int^t_0\nabla b(X_s) Y_s\dif s+A_1 \int^t_0h(s)\dif s+ A_2\int^t_0\!\!\!\int_{\Gamma_0}\v(s,z)N(\dif s,\dif z).
$$
By Fatou's lemma and (\ref{EW10}), we have
$$
\varlimsup_{n\to\infty}\mE|D_\Theta X^n_t-Y_t|^p\leq \|\nabla b\|^p_\infty\int^t_0\varlimsup_{n\to\infty}\mE|D_\Theta X^{n-1}_s-Y_s|^p\dif s,
$$
which then gives
$$
\varlimsup_{n\to\infty}\mE|D_\Theta X^n_t-Y_t|^p=0.
$$
Thus, by (\ref{EW10}) we have $X_t\in\mW^{1,p}_\Theta$ and $D_\Theta X_t=Y_t$.
The proof is complete.
\end{proof}
Let $J_t:=J_t(x):=\nabla X_t(x)$ be the Jacobii's matrix and $K_t:=K_t(x):=J^{-1}_t(x)$. Then $J_t$ and $K_t$ solve the following ODEs
\begin{align}
J_t=\mI+\int^t_0\nabla b(X_s) J_s\dif s,\ \ K_t=\mI-\int^t_0K_s\nabla b(X_s)\dif s,\label{EY11}
\end{align}
and it is easy to see that
\begin{align}
\sup_{t\in[0,1]}\sup_{x\in\mR^d}|J_t(x)|\leq\e^{\|\nabla b\|_\infty},\ \ \sup_{t\in[0,1]}\sup_{x\in\mR^d}|K_t(x)|\leq\e^{\|\nabla b\|_\infty}.\label{TY1}
\end{align}
By (\ref{EW6}) and the formula of constant variation, we have
\begin{align}
D_\Theta X_t=J_t\int^t_0K_sA_1 h(s)\dif s+J_t\int^t_0\!\!\!\int_{\Gamma_0}K_sA_2\v(s,z) N(\dif s,\dif z).\label{EW7}
\end{align}
Below, let $\zeta(z)$ be a nonnegative smooth function with
\begin{align}
\zeta(z)=|z|^3,\ \ |z|\leq 1/4,\ \ \zeta(z)=0,\ \ |z|>1/2.\label{Ze}
\end{align}
Let us choose
$$
\Theta_j(x)=(h_j(x;\cdot),\v_j(x;\cdot))
$$
with
$$
h_j(x; s)=(K_s(x)A_1)^*_{\cdot j}, \ \ \v_j(x;s,z)=(K_s(x)A_2)^*_{\cdot j}\zeta(z).
$$
\bl
Under {\bf (H$^\alpha_1$)}, for each $j=1,\cdots,d$ and $x\in\mR^d$, we have $\Theta_j(x)\in\mH_{\infty-}\times\mV_{\infty-}$ and
\begin{align}
\div\Theta_j(x)=-\sum_l\int^1_0(K_s(x)A_1)_{lj}\dif W^l_s+\sum_l\int^1_0\!\!\!\int_{\Gamma_0}(K_s(x)A_2)_{lj}\eta_l(z)\tilde N(\dif z,\dif s),\label{TY88}
\end{align}
where $\eta_l(z):=\p_l\zeta(z)+\zeta(z)\p_l\log\kappa(z)$. In particular, for any $p\geq 2$,
\begin{align}
\sup_{x\in\mR^d}\mE|\div\Theta_j(x)|^p<+\infty.\label{213}
\end{align}
\el
\begin{proof}
Since $\dis(z,\Gamma_0^c)\geq |z|\wedge(1-|z|)$,  by (\ref{TY1}) and (\ref{Ze}), it is easy to check that $\Theta_j(x)\in\mH_{\infty-}\times\mV_{\infty-}$.
Moreover,  by definition (\ref{ER1}) we immediately have (\ref{TY88}). As for (\ref{213}), it follows by (\ref{214}).
\end{proof}
Write
$$
{\bf\Theta}:=(\Theta_1,\cdots,\Theta_d),\ \ (D_{\bf\Theta} X_t)_{ij}:=D_{\Theta_j}X^i_t
$$
and
\begin{align}
\Sigma_t(x):=\int^t_0K_s(x)A_1A_1^*K_s^*(x)\dif s+\int^t_0\!\!\!\int_{\Gamma_0}K_s(x) A_2A^*_2K_s^*(x)\zeta(z) N(\dif s,\dif z),\label{TY7}
\end{align}
then by (\ref{EW7}),
\begin{align}
D_{\bf\Theta} X_t(x)=J_t(x)\Sigma_t(x).\label{TY8}
\end{align}
The matrix $\Sigma_t(x)$ will be called the reduced Malliavin matrix (cf. \cite[p. 89, (7-20)]{Bi-Gr-Ja} and \cite[(2.12)]{Xu}).
\bl\label{Le32}
Assume that $b\in C^\infty$ has bounded derivatives of all orders.
For any $k,n\in\mN\cup\{0\}$ with $k+n\geq 1$, $j_1,\cdots, j_n\in\{1,\cdots,d\}$ and $p\geq 2$, we have
\begin{align}
&\sup_{t\in[0,1]}\sup_{x\in\mR^d}\mE|D_{\Theta_{j_1}}\cdots D_{\Theta_{j_n}}\nabla^k X_t(x)|^p<\infty,\label{TY3}\\
&\sup_{t\in[0,1]}\sup_{x\in\mR^d}\mE|D_{\Theta_{j_1}}\cdots D_{\Theta_{j_n}}J_t(x)|^p<\infty,\label{TY4}\\
&\sup_{t\in[0,1]}\sup_{x\in\mR^d}\mE|D_{\Theta_{j_1}}\cdots D_{\Theta_{j_n}}K_t(x)|^p<\infty,\label{TY5}\\
&\sup_{t\in[0,1]}\sup_{x\in\mR^d}\mE|D_{\Theta_{j_1}}\cdots D_{\Theta_{j_n}}\Sigma_t(x)|^p<\infty.\label{TY6}
\end{align}
Moreover, under {\bf (H$^\alpha_m$)} with $m\geq 2$, for any $n\leq m-1$, we also have
\begin{align}
\sup_{t\in[0,1]}\sup_{x\in\mR^d}\mE|D_{\Theta_{j_1}}\cdots D_{\Theta_{j_n}}\div\Theta_i(x)|^p<\infty.\label{TY77}
\end{align}
\el
\begin{proof}
First of all, by equation (\ref{EY11}) and the induction, it is easy to prove that for any $k\in\mN$,
\begin{align}
\sup_{t\in[0,1]}\sup_{x\in\mR^d}\sup_{\omega}|\nabla^k X_t(x,\omega)|<+\infty.\label{TY66}
\end{align}
By (\ref{TY7}), (\ref{TY1}) and  inequality (\ref{BT1}), we have for any $p\geq 1$,
$$
\sup_{t\in[0,1]}\sup_{x\in\mR^d}\mE|\Sigma_t(x)|^p<\infty,
$$
which together with (\ref{TY8}) and (\ref{TY1}) yields (\ref{TY3}) with $n=1$ and $k=0$. By induction,
the higher order derivatives for (\ref{TY3})-(\ref{TY6})  follow by (\ref{EW6}), (\ref{TY7}), Proposition \ref{Pr1} and (\ref{BT1}).

We now look at (\ref{TY77}). By (\ref{TY88}) and Proposition \ref{Pr1}, we have
\begin{align*}
D_{\Theta_j}\div\Theta_i&=-\sum_l\int^1_0(D_{\Theta_j}K_sA_1)_{li}\dif W^l_s-\sum_l\int^1_0(K_sA_1)_{li}(K_sA_1)_{lj}\dif s\\
&\quad+\sum_{l,l'}\int^1_0\!\!\!\int_{\Gamma_0}(K_sA_2)_{li}(K_sA_2)_{l'j}\p_{l'}\eta_l(z)\zeta(z)N(\dif z,\dif s)\\
&\quad+\sum_l\int^1_0\!\!\!\int_{\Gamma_0}(D_{\Theta_j}K_sA_2)_{li}\eta_l(z)\tilde N(\dif z,\dif s).
\end{align*}
Recalling $\eta_l(z):=\p_l\zeta(z)+\zeta(z)\p_l\log\kappa(z)$ and $\zeta(z)$ given by (\ref{Ze}), by {\bf (H$^\alpha_m$)} we have
$$
|\eta_l(z)|\leq C|z|^2,\ \ |\p_{l'}\eta_l(z)\zeta(z)|\leq C|z|^4,\ \ z\in\Gamma_0.
$$
By (\ref{BT1}) and (\ref{TY5}), we obtain (\ref{TY77}) for $n=1$. The higher order derivatives follow by induction.
\end{proof}
\section{Proof of Theorem \ref{Main1}}
\subsection{Invertibility of $\Sigma_t$}
We first prove the following lemma as in  \cite[Lemma 2.1]{Zh2}.
\bl\label{Le51}
Set $\Delta L^0_s:=L^0_s-L^0_{s-}$ and define
$$
\Omega_0:=\Big\{\omega: \{s: |\Delta L^0_s(\omega)|\not=0\} \mbox{ is dense in $[0,\infty)$}\Big\}.
$$
Under {\bf (H$^\alpha_1$)}, we have $\mP(\Omega_0)=1.$
\el
\begin{proof}
Define a stopping time $\tau:=\inf\{t>0: |L^0_t|=0\}$. As in the proof of \cite[Lemma 2.1]{Zh2}, it suffices to prove that
$$
\mP(\tau=0)=1.
$$
For any $\eps\in(0,1)$, we have
\begin{align*}
\eps^2\geq \mE\left(|\Delta L^0_\tau|^21_{|\Delta L^0_\tau|\leq\eps}\right)
&=\mE\left(\sum_{0<s\leq\tau}|\Delta L^0_s|^21_{|\Delta L^0_s|\leq\eps}\right)=\mE\left(\int^{\tau}_0\!\!\!\int_{|z|\leq\eps}|z|^2N(\dif s,\dif z)\right)\\
&=\mE\left(\int^{\tau}_0\!\!\!\int_{|z|\leq\eps}|z|^2\nu(\dif z)\dif s\right)
=\int_{|z|\leq\eps}|z|^2\kappa(z)\dif z\mE\tau,
\end{align*}
which, together with {\bf (H$^\alpha_\kappa$)} and letting $\eps\to 0$, implies that
$$
\mE\tau=0\Rightarrow \mP(\tau=0)=1.
$$
The proof is complete.
\end{proof}
Let $V:\mR^d\to\mR^d\times\mR^d$ be a matrix-valued $C^2_b$-function. Below we set
$$
M^V_t:=\sum_{l=1}^d\int^t_0K_s\p_l V(X_{s-})\dif W^l_s,\  \ G^V_t:=\int^t_0\!\!\!\int_{\Gamma_0}K_s(A_2z\cdot \nabla) V(X_{s-})\tilde N(\dif s,\dif z),
$$
and
$$
H^V_t:=\sum_{0<s\leq t}K_s\Big(V(X_s)-V(X_{s-})-(\Delta X_s\cdot\nabla) V(X_{s-})\Big).
$$
We have
\bl
There exists a subsequence $n_m\to\infty$ such that $\mP(\Omega^V_1)=1$, where
\begin{align*}
\Omega^V_1:=\Bigg\{\omega: &H^V_t=\lim_{n_m\to\infty}\sum_{0<s\leq t}
K_s\Big(V(X_s)-V(X_{s-})-(\Delta X_s\cdot\nabla V)(X_{s-}) \Big)1_{|\Delta L^0_s|>\frac{1}{n_m}}\mbox{ and }\\
&G^V_t=\lim_{n_m\to\infty}\sum_{s\in(0,t]}K_s(A_2\Delta L^0_s\cdot \nabla) V(X_{s-})
1_{|\Delta L^0_s|>\frac{1}{n_m}}\mbox{ uniformly in $t\in[0,1]$.}\Bigg\}.
\end{align*}
\el
\begin{proof}
By $\nu(\dif z)=\nu(-\dif z)$ and Doob's maximal inequality, we have
\begin{align*}
&\mE\left(\sup_{t\in[0,1]}\left|G^V_t-\sum_{s\in(0,t]}K_s(A_2\Delta L^0_s\cdot \nabla) V(X_{s-})1_{|\Delta L^0_s|>\frac{1}{n}}\right|^2\right)\\
&\qquad=\mE\left(\sup_{t\in[0,1]}\left|\int^t_0\!\!\!\int_{|z|\leq\frac{1}{n}}K_s(A_2z\cdot \nabla) V(X_{s-})\tilde N(\dif s,\dif z)\right|^2\right)\\
&\qquad\leq 4\mE\left|\int^1_0\!\!\!\int_{|z|\leq\frac{1}{n}}K_s(A_2z\cdot \nabla) V(X_{s-})\tilde N(\dif s,\dif z)\right|^2\\
&\qquad\leq 4\mE\left(\int^1_0\!\!\!\int_{|z|\leq\frac{1}{n}}|K_s(A_2z\cdot \nabla) V(X_{s-})|^2\nu(\dif z)\dif s\right)\\
&\qquad\leq C\int_{|z|\leq\frac{1}{n}}|z|^2\nu(\dif z)\to 0,\ \ n\to\infty.
\end{align*}
Similarly, we also have
\begin{align*}
&\mE\left(\sup_{t\in[0,1]}\left|\sum_{0<s\leq t}K_s\Big(V(X_s)-V(X_{s-})-(\Delta X_s\cdot \nabla) V(X_{s-})\Big)1_{|\Delta L^0_s|\leq\frac{1}{n}}\right|\right)\\
&\qquad=\mE\left(\sup_{t\in[0,1]}\left|\int^t_0\!\!\!\int_{|z|\leq\frac{1}{n}}K_s\Big(V(X_{s-}+z)-V(X_{s-})-(A_2z\cdot\nabla V)(X_{s-})\Big)N(\dif s,\dif z)\right|\right)\\
&\qquad\leq\mE\left(\int^1_0\!\!\!\int_{|z|\leq\frac{1}{n}}|K_s(V(X_{s-}+z)-V(X_{s-})-(A_2z\cdot\nabla) V(X_{s-}))|\nu(\dif z)\dif s\right)\\
&\qquad\leq C\int_{|z|\leq\frac{1}{n}}|z|^2\nu(\dif z)\to 0,\ \ n\to\infty.
\end{align*}
The proof is complete.
\end{proof}
By \cite[p.64, Theorem 21 and p.68, Theorem 23]{Pr}, we have
\bl
For $n\in\mN$, let $t_k:=k/n$. There exists a subsequence $n_m\to\infty$ such that 
$$
\mP(\Omega^V_2\cap\Omega^V_3)=1,
$$ 
where
\begin{align*}
\Omega^V_2&:=\Bigg\{\omega:
\lim_{n_m\to\infty}\sum_{k=0}^{n_m-1}\sum_lK_{t_k}\p_l V(X_{t_k})(W^l_{t_{k+1}\wedge t}-W^l_{t_k\wedge t})=M^V_t\mbox{ uniformly in $t\in[0,1]$}\Bigg\},\\
\Omega^V_3&:=\Bigg\{\omega:\  \forall i,j,i',j'=1,\cdots,d, 
\lim_{n_m\to\infty}\sum_{k=0}^{n_m-1}(M^V_{t_{k+1}\wedge t}-M^V_{t_k\wedge t})_{ij}(M^V_{t_{k+1}\wedge t}-M^V_{t_k\wedge t})_{i'j'}\\
&\qquad\quad=\sum_{k,l,k'}\int^t_0(K_s)_{ik}\p_lV_{jk}(X_s)(K_s)_{i'k'}\p_lV_{j'k'}(X_s)\dif s\mbox{ uniformly in $t\in[0,1]$}\Bigg\}.
\end{align*}
\el
By It\^o's formula, we also have
\bl\label{Le54}
Let $[b,V]:=b\cdot\nabla V-\nabla b\cdot V+\tfrac{1}{2}\nabla^2_{A_1A_1^*} V$ and define
\begin{align*}
\Omega^V_4:=\Bigg\{&\omega: K_tV(X_t)=V(x)+\int^t_0K_s[b,V](X_s)\dif s+H^V_t+M^V_t+G^V_t,\forall t\geq 0\Bigg\}.
\end{align*}
Then $\mP(\Omega^V_4)=1$.
\el
Now we can prove the following key lemma.
\bl\label{Le45}
Fix $x\in\mR^d$. Let $B_0:=\mI$ and for $n\in\mN$,
$$
B_n(x):=[b,B_{n-1}](x):=b(x)\cdot\nabla B_{n-1}(x)-\nabla b(x)\cdot B_{n-1}(x)+\tfrac{1}{2}\nabla^2_{A_1A_1^*}B_{n-1}(x).
$$ 
Assume that for some $n=n(x)\in\mN$,
\begin{align}
\mathrm{Rank}[A_1, B_1(x)A_1,\cdots, B_n(x)A_1, A_2, B_1(x)A_2,\cdots, B_n(x)A_2]=d.\label{QW1}
\end{align}
Then under {\bf (H$^\alpha_1$)},  for any $t>0$,  $\Sigma_t(x)$ is almost surely invertible.
\el
\begin{proof}
Set
$$
\tilde\Omega:=\cap_{n=1}^\infty(\Omega^{B_n}_1\cap\Omega^{B_n}_2\cap\Omega^{B_n}_3\cap\Omega^{B_n}_4)\cap\Omega_0.
$$
Then by Lemmas \ref{Le51}-\ref{Le54}, we have
$$
\mP(\tilde\Omega)=1.
$$
We want to prove that under (\ref{QW1}), for each $t>0$, the reduced Malliavin matrix $\Sigma_t(x,\omega)$ is invertible for each $\omega\in\tilde\Omega$.
Without loss of generality, we assume $t=1$ and fix an $\omega\in\tilde\Omega$. For simplicity of notation, we shall drop $(x,\omega)$ below.
By (\ref{TY7}), for a row vector $u\in\mR^d$ we have
\begin{align*}
u\Sigma_1 u^*&=\int^1_0|uK_sA_1|^2\dif s+\int^1_0\!\!\!\int_{\Gamma_0}|uK_sA_2|^2\zeta(z) N(\dif s,\dif z)\\
&=\int^1_0|uK_sA_1|^2\dif s+\sum_{s\leq 1}|uK_sA_2|^2\zeta(\Delta L^0_s) 1_{|\Delta L^0_s|\not=0}.
\end{align*}
Suppose that for some $u\in\mS^{d-1}$,
$$
u\Sigma_1 u^*=0.
$$
Since $s\mapsto K_s$ is continuous and $\omega\in\Omega_0$, we have
$$
|uK_sA_1|^2=|uK_s A_2|^2=0,\ \ \forall s\in[0,1].
$$
Hence, by (\ref{EY11}) we have
$$
0=uK_tA_i=uA_i-\int^t_0uK_s\nabla b(X_s)A_i\dif s,\forall t\in[0,1], \ i=1,2,
$$
which implies that
$$
uA_i=0,\ \ i=1,2,
$$
and
$$
uK_t\nabla b(X_t)A_i=uK_tB_1(X_t)A_i=0,\ \ t\in[0,1],\ \ i=1,2.
$$
Now we use the induction to prove
\begin{align}
uK_tB_n(X_t)A_i=0,\ \ t\in[0,1],\ \ i=1,2.\label{In}
\end{align}
Suppose that (\ref{In}) holds for some $n\in\mN$.
In view of $\omega\in\Omega^{B_n}_4$, we have for all $t\in[0,1]$,
\begin{align*}
uK_tB_n(X_t)A_i&=uB_n(x)A_i+\int^t_0uK_sB_{n+1}(X_s)A_i\dif s+uH^{B_n}_tA_i+uM^{B_n}_tA_i+uG^{B_n}_tA_i.
\end{align*}
By the induction hypothesis and the definition of $H^{B_n}_t$, we further have
\begin{align}
0&=\int^t_0uK_sB_{n+1}A_i(X_s)\dif s-\sum_{0<s\leq t}uK_s(\Delta X_s\cdot\nabla) B_n(X_{s-}) A_i+uM^{B_n}_tA_i+uG^{B_n}_tA_i\no\\
&=\int^t_0uK_sB_{n+1}A_i(X_s)\dif s+uM^{B_n}_tA_i,\ \forall t\in[0,1],\quad (\because \omega\in\Omega^{B_n}_1),\label{QW2}
\end{align}
which together with $\omega\in\Omega^{B_n}_3$ implies that
$$
0=\lim_{n_m\to\infty}\sum_{k=0}^{n_m-1}\<uM^{B_n}_{t_{k+1}}A_i-uM^{B_n}_{t_{k}}A_i, uM^{B_n}_{t_{k+1}}A_i-uM^{B_n}_{t_{k}}A_i\>_{\mR^d}
=\sum_l\int^1_0|uK_s\p_l B_n(X_s)A_i|^2\dif s.
$$
In particular,
$$
uK_s\p_l B_n(X_s)A_i=0,\ \ \forall s\in[0,1].
$$
Since $\omega\in\Omega^{B_n}_2$, we also have
$$
uM^{B_n}_tA_i=0,\ \  \forall t\in[0,1],
$$
which together with (\ref{QW2}) implies that
$$
uK_sB_{n+1}A_i(X_s)=0,\ \ \forall s\in[0,1].
$$
Thus, we obtain
$$
uA_i=uB_1A_i=\cdots uB_nA_i=0,\ \ i=1,2,
$$
which is contradict with (\ref{QW1}). The proof is complete.
\end{proof}
\subsection{Proof of Theorem \ref{Main1}}
Now we can finish the proof of Theorem \ref{Main1} by the same argument as in \cite{Do-Pe-So-Zh}.
We divide the proof into two steps.
\\
\\
(1) Let GL$(d)\simeq\mR^d\times\mR^d$ be the set of all $d\times d$-matrix. Define
$$
\mM_n:=\Big\{\Sigma\in \mathrm{GL}(d): |\Sigma|\leq n,\ \ \det(\Sigma)\geq1/n\Big\}.
$$
Then $\mM_n$ is a compact subset of GL$(d)$. Let $\Phi_n\in C^\infty(\mR^d\times\mR^d)$ be a smooth function so that
$$
\Phi_n|_{\mM_n}=1,\ \ \Phi_n|_{\mM^c_{n+1}}=0,\ \ 0\leq \Phi_n\leq 1.
$$
Below we fix $t>0$ and $x\in\mR^d$. For each $n\in\mN$, let us define a finite measure $\mu_n$ by
$$
\mu_n(A):=\mE\Big[1_A(X_t)\Phi_n(\Sigma_t )\Big],\ \ A\in\sB(\mR^d).
$$
For each $\varphi\in C^\infty_b(\mR^d)$, by the chain rule and (\ref{TY8}), we have
$$
D_{\bf \Theta}(\varphi(X_t))=\nabla\varphi(X_t)D_{\bf\Theta}X_t=\nabla\varphi(X_t)J_t\Sigma_t,
$$
where $\nabla=(\p_1,\cdots,\p_d)$. So,
\begin{align}
\nabla\varphi(X_t)=D_{\bf \Theta}(\varphi(X_t))\Sigma^{-1}_tK_t.\label{WQ2}
\end{align}
Thus, by the integration by parts (\ref{ER88}), we have for $i=1,\cdots,d$,
\begin{align}
\int_{\mR^d}\p_i\varphi(y)\mu_n(\dif y)&=\mE[\p_i\varphi(X_t)\Phi_n(\Sigma_t)]\no\\
&=\sum_j\mE\Big[D_{\Theta_j}(\varphi(X_t))(\Sigma^{-1}_tK_t)_{ij}\Phi_n(\Sigma_t)\Big]=\mE[\varphi(X_t) H^i_t],\label{WQ3}
\end{align}
where
$$
H^i_t:=-\sum_j\left((\Sigma^{-1}_tK_t)_{ij}\Phi_n(\Sigma_t)\div(\Theta_j)+D_{\Theta_j}((\Sigma^{-1}_tK_t)_{ij}\Phi_n(\Sigma_t))\right).
$$
From this and using Lemma \ref{Le32}, by cumbersome calculations, we derive that
$$
\left|\int_{\mR^d}\p_i\varphi(y)\mu_n(\dif y)\right|\leq C_n\|\varphi\|_\infty,\ \ i=1,\cdots,d,
$$
where $C_n$ is independent of $t,x$. Hence, $\mu_n$ is absolutely continuous with respect to the Lebesgue measure
(cf. \cite{Nu}), and by the Sobolev embedding theorem (cf. \cite{Ad}), the density $p_n(y)$ satisfies that for any $q\in[1,d/(d-1))$,
$$
\int_{\mR^d}p_n(y)^q\dif y\leq C_{d,q,n},
$$
where the constant $C_{d,q,n}$ is independent of $t,x$. Therefore, for any Borel set $F\subset\mR^d$ and $R>0$,  we have
\begin{align}
\mu_n(F)=\int_Fp_n(y)\dif y&\leq m(F)R+\int_{F\cap \{p_n>R\}}p_n(y)\dif y\leq m(F)R+\frac{C_{d,q,n}}{R^{q-1}},\label{ET4}
\end{align}
where $m$ is the Lebesgue measure and $q>1$. In particular, for Lebesgue zero measure $A\subset\mR^d$, 
$$
\mE\Big[1_A(X_t)\Phi_n(\Sigma_t )\Big]=0.
$$
By Lemma \ref{Le45} and the dominated convergence therem, we obtain that for any Lebesgue zero measure $A\subset\mR^d$, 
$$
\mE[1_A(X_t)]=0,
$$
which means that the law of $X_t$ is absolutely continuous with respect to the Lebesgue measure.
\\
\\
(2) Let $\chi_n\in C^\infty(\mR^d)$ be a smooth function with
\begin{align}
\chi_n|_{\{|x|\leq n\}}=1,\ \ \chi_n|_{\{|x|>n+1\}}=0, \ \ 0\leq\chi_n\leq 1.\label{Chi}
\end{align}
Let $f$ be a bounded nonnegative measurable function. By Lusin's theorem, for any $\eps>0$, 
there exist a set $F_\eps\subset \{x\in\mR^d:|x|<n+1\}$ and a nonnegative continuous function $g\in C_c(\mR^d)$ such that
$$
f\chi_n|_{F^c_\eps}=g|_{F^c_\eps}, \ \ \|g\|_\infty\leq \|f\|_\infty,\ \ m(F_\eps)<\eps.
$$
Let $\mu_{t,x;n}$ be defined by
$$
\mu_{t,x;n}(A):=\mE\Big[1_A(X_t(x))\Phi_n(\Sigma_t(x))\Big],\ \ A\in\sB(\mR^d).
$$
By the dominated convergence theorem and (\ref{ET4}), we have for any $R>0$,
\begin{align*}
&\varlimsup_{x\to x_0}\mE\Big[(f\chi_n)(X_t(x))\Phi_n(\Sigma_t(x))\Big]\\
&\quad\leq\varlimsup_{x\to x_0}\mE\Big[g(X_t(x))\Phi_n(\Sigma_t(x))\Big]
+\varlimsup_{x\to x_0}\mE\Big[|f\chi_n-g|(X_t(x))\Phi_n(\Sigma_t(x))\Big]\\
&\quad\leq \mE\Big[g(X_t(x_0))\Phi_n(\Sigma_t(x_0))\Big]+2\|f\|_\infty\varlimsup_{x\to x_0}\mu_{t,x;n}(F_\eps)\\
&\quad\leq\mE\Big[(f\chi_n)(X_t(x_0))\Phi_n(\Sigma_t(x_0))\Big]+\mE\Big[(g-f\chi_n)(X_t(x_0))\Phi_n(\Sigma_t(x_0))\Big]\\
&\qquad+2\|f\|_\infty\left(m(F_\eps)R+\frac{C_{d,q,n}}{R^{q-1}}\right)\\
&\quad\leq \mE f(X_t(x_0))+4\|f\|_\infty\left(m(F_\eps)R+\frac{C_{d,q,n}}{R^{q-1}}\right).
\end{align*}
First letting $\eps\to 0$ and then $R\to\infty$, we obtain for $n\in\mN$,
\begin{align}
\varlimsup_{x\to x_0}\mE\Big[(f\chi_n)(X_t(x))\Phi_n(\Sigma_t(x))\Big]\leq \mE f(X_t(x_0)).\label{ET3}
\end{align}
On the other hand, by the definition (\ref{TY7}) of $\Sigma_t(x)$, it is easy to see that
$$
x\mapsto X_t(x),\ \Sigma_t(x) \mbox{ are continuous in probability}.
$$ 
Thus, by the dominated convergence theorem and (\ref{ET3}), we have
\begin{align}
\varlimsup_{x\to x_0}\mE f(X_t(x))&\leq \varlimsup_{x\to x_0}\mE\Big[(f\chi_n)(X_t(x) )\Phi_n(\Sigma_t(x) )\Big]
+\|f\|_\infty\varlimsup_{x\to x_0}\mE\Big[1-\chi_n(X_t(x) )\Phi_n(\Sigma_t(x) )\Big]\no\\
&= \varlimsup_{x\to x_0}\mE\Big[(f\chi_n)(X_t(x) )\Phi_n(\Sigma_t(x) )\Big]
+\|f\|_\infty\mE\Big[1-\chi_n(X_t(x_0) )\Phi_n(\Sigma_t(x_0))\Big]\no\\
&\leq \mE f(X_t(x_0))+\|f\|_\infty\mP\Big(\{\Sigma_t(x_0)\notin \mM_n\}\cup\{|X_t(x_0)|>n\}\Big)\label{ET2},
\end{align}
which, by Lemma \ref{Le45} and letting $n\to\infty$, implies
$$
\varlimsup_{x\to x_0}\mE f(X_t(x))\leq\mE f(X_t(x_0)).
$$
Applying the above limit to the nonnegative function $\|f\|_\infty-f(x)$, we also have
$$
\varlimsup_{x\to x_0}\mE (\|f\|_\infty-f(X_t(x)))\leq\|f\|_\infty-\mE f(X_t(x_0))\Rightarrow
\varliminf_{x\to x_0}\mE f(X_t(x))\geq\mE f(X_t(x_0)).
$$
Thus, we obtain the desired continuity (\ref{WQ1}).
\section{Proof of Theorem \ref{Main2}}

\subsection{Norris' type estimate}
We first recall the following Norris' type estimate (cf. \cite{Zh3}).
\bl\label{Lemma2}
Let $Y_t=y+\int^t_0\beta_s\dif s$ be an $\mR^d$-valued process, where $\beta_t$ takes the following form:
$$
\beta_t=\beta_0+\int^t_0\gamma_s\dif s+\int^t_0 Q_s\dif W_s+\int^t_0\!\!\!\int_{\Gamma_0}g_s(z)\tilde N(\dif s,\dif z),
$$
where $\gamma_t:\mR_+\to\mR^d$, $ Q_t:\mR_+\to\mR^d\times\mR^d$ and $g_t(z):\mR_+\times\mR^d\to\mR^d$ are three left continuous $\sF_t$-adapted processes.
Suppose that for some $R>0$,
\begin{align}
|\beta_t|, |Q_t|, |\gamma_t|\leq R,\ \ |g_t(z)|\leq R(1\wedge|z|).
\end{align}
Then there exists a constant $C\geq1$ such that for any $t\in(0,1)$, $\delta\in(0,\frac{1}{3})$ and $\eps\in(0,t^3)$,
\begin{align}
P\left\{\int^t_0|Y_s|^2\dif s<\eps, \int^t_0|\beta_s|^2\dif s\geq 9R^2\eps^{\delta}\right\}\leq
4\exp\left\{-\frac{\eps^{\delta-\frac{1}{3}}}{CR^4}\right\}.\label{EW5}
\end{align}
\el
The following lemma is simple.
\bl
Assume that for some $\alpha\in(0,2)$,
\begin{align}
\lim_{\eps\to 0}\eps^{\alpha-2}\int_{|z|\leq\eps}|z|^2\nu(\dif z)=:c_1>0.\label{UY21}
\end{align}
Then for any $p\geq 2$, there exist constants $\eps_0, c_2>0$ such that for all $\eps\in(0,\eps_0)$,
\begin{align}
\int_{|z|\leq\eps}|z|^p\nu(\dif z)\geq c_2\eps^{p-\alpha}.\label{UY2}
\end{align}
\el
\begin{proof}
For any $\delta\in(0,1)$, by (\ref{UY21}), there is an $\eps_0>0$ such that for all $\eps\in(0,\eps_0)$,
$$
(1-\delta)c_1\eps^{2-\alpha}\leq\int_{|z|\leq\eps}|z|^2\nu(\dif z)\leq (1+\delta)c_1\eps^{2-\alpha}.
$$
Hence,
\begin{align*}
\int_{|z|\leq\eps}|z|^p\nu(\dif z)&=\sum_{n=0}^\infty\int_{2^{-(n+1)}\eps<|z|\leq 2^{-n}\eps}|z|^p\nu(\dif z)
\geq\sum_{n=0}^\infty (2^{-(n+1)}\eps)^{p-2}\int_{2^{-(n+1)}\eps<|z|\leq 2^{-n}\eps}|z|^2\nu(\dif z)\\
&\geq\sum_{n=0}^\infty (2^{-(n+1)}\eps)^{p-2}\left((1-\delta)c_1(2^{-n}\eps)^{2-\alpha}-(1+\delta)c_1(2^{-(n+1)}\eps)^{2-\alpha}\right)\\
&=\eps^{p-\alpha}c_12^{\alpha-p}\sum_{n=0}^\infty 2^{-n(p-\alpha)}\left(2^{2-\alpha}(1-\delta)-(1+\delta)\right),
\end{align*}
which gives (\ref{UY2}) by letting $\delta$ small enough.
\end{proof}
We also need the following estimate.
\bl\label{Le52}
Let $g_t$ be a nonnegative bounded predictable processes. Under {\bf (H$^\alpha_1$)}, 
there exist constants $\lambda_0, c_0\geq 1$ depending on the bound of $g_t$ such that for all $\lambda>\lambda_0$,
\begin{align}
\mP\left\{\int^t_0\!\!\!\int_{\Gamma_0}g_s\zeta(z)N(\dif s,\dif z)\leq\eps; \int^t_0g_s\dif s>\eps^{\frac{\alpha}{6}}\right\}
\leq \exp\left\{1-c_0\eps^{-\frac{\alpha}{6}}\right\},
\end{align}
where $\zeta(z)$ is defined by (\ref{Ze}).
\el
\begin{proof}
Define
$$
\beta^\lambda_t:=\int_{\Gamma_0}(1-\e^{-\lambda g_t\zeta(z)})\nu(\dif z)
$$
and
$$
M^\lambda_t:=-\lambda\int^t_0\!\!\!\int_{\Gamma_0}g_s\zeta(z)N(\dif s,\dif z)+\int^t_0\beta^\lambda_s\dif s.
$$
By It\^o's formula, we have
$$
\e^{M^\lambda_t}=1+\int^t_0\!\!\!\int_{\Gamma_0}\e^{M^\lambda_{s-}}(\e^{-\lambda g_s\zeta(z)}-1)\tilde N(\dif s,\dif z).
$$
Since for any $x\geq 0$,
$$
1-\mathrm{e}^{-x}\leq 1\wedge x,
$$
we have
$$
M^\lambda_t\leq \int^t_0\beta^\lambda_s\dif s
\leq \int^t_0\!\!\!\int_{\Gamma_0}(1\wedge(\lambda g_s\zeta(z)))\nu(\dif z)\dif s.
$$
Hence,
\begin{align}
\mE\e^{M^\lambda_t}=1.\label{ERP3}
\end{align}
Since for any $\kappa\in(0,1)$ and $x\leq-\log k$,
$$
1-\mathrm{e}^{-x}\geq \kappa x,
$$
letting $\kappa=\frac{1}{\e}$ and by (\ref{UY2}), there exist  $\lambda_0,c_0>1$ such that for all $\lambda\geq\lambda_0$,
\begin{align*}
\beta^{\lambda}_s\geq\int_{|z|\leq((\|g\|_\infty+1)\lambda)^{-1/3}}(1-\e^{-\lambda g_s\zeta(z)})\nu(\dif z)
\geq\frac{\lambda g_s}{\e}\int_{|z|\leq((\|g\|_\infty+1)\lambda)^{-1/3}}|z|^3\nu(\dif z)\geq c_0\lambda^{\frac{\alpha}{3}} g_s.
\end{align*}
Thus, 
\begin{align*}
&\left\{\int^t_0\!\!\!\int_{\Gamma_0}g_s\zeta(z)N(\dif s,\dif z)\leq\eps; \int^t_0g_s\dif s>\eps^{\frac{\alpha}{6}}\right\}\\
&\quad=\left\{\e^{M^\lambda_t}\geq \e^{-\lambda\eps+\int^t_0\beta^\lambda_s\dif s}; \int^t_0g_s\dif s>\eps^{\frac{\alpha}{6}}\right\}
\subset\left\{\e^{M^\lambda_t}\geq \e^{-\lambda\eps+c_0\lambda^{\frac{\alpha}{3}}\eps^{\frac{\alpha}{6}}}\right\},
\end{align*}
which, by Chebyschev's inequality, (\ref{ERP3}) and letting $\lambda=\frac{1}{\eps}$, gives the desired estimate.
\end{proof}
\subsection{Proof of Theorem \ref{Main2}}
\bl\label{Le53}
Under (\ref{Con}), there exists constants $C_1, C_2\in(0,1)$ independent of the starting point $x$ and $t_0\in(0,1)$ such that for all $t\in(0,t_0)$ and $\eps\in(0, C_1t^4)$,
\begin{align}
\sup_{|u|=1}\mP\left(\int^t_0(|uK_sA_1|^2+|uK_sA_2|^2)\dif s\leq\eps\right)\leq 8\exp\left\{-C_2\eps^{-\frac{1}{12}}\right\}.\label{EW3}
\end{align}
\el
\begin{proof}
Fix $u\in\mS^{d-1}$ and set for $i=1,2$ and $j=1,\cdots, d$,
\begin{align*}
&Y^i_t:=uK_tA_i,\quad\beta^i_t:=uK_t\nabla b(X_t)A_i,\quad  Q^{ij}_t:=\sum_kuK_t(A_1)_{kj}\p_j\nabla b(X_t)A_i, \\
&\gamma^i_t:=uK_t\Bigg[\left((b\cdot\nabla)\nabla b-(\nabla b)^2+\tfrac{1}{2}\nabla^2_{A_1A_1^*}\nabla b\right)(X_t)\\
&\qquad+\int_{\Gamma_0}\Big(\nabla b(X_t+A_2z)-\nabla b(X_t)-1_{|z|\leq 1}(A_2z\cdot\nabla)\nabla b(X_t)\Big)\nu(\dif z)\Bigg]A_i,\\
&g^i_t(z):=uK_t(\nabla b(X_{t-}+A_2z)-\nabla b(X_{t-}))A_i.
\end{align*}
By equations (\ref{EY11}) and It\^o's formula, one sees that 
$$
Y^i_t=uA_i+\int^t_0\beta^i_s\dif s,
$$ 
and
$$
\beta^i_t=u\nabla b(x)A_i+\int^t_0\gamma^i_s\dif s+\sum_j\int^t_0 Q^{ij}_s\dif W^j_s+\int^t_0\!\!\!\int_{\Gamma_0}g^i_s(y)\tilde N(\dif s,\dif y).
$$
By the assumptions, it is easy to see that for some $R>0$,
$$
|g^i_t(z)|\leq R(1\wedge|z|),\ \ |\beta^i_t|+|\gamma^i_t|+|Q^{ij}_t|\leq R.
$$
Notice that
\begin{align*}
&\left\{\int^t_0(|Y^1_s|^2+|Y^2_s|^2)\dif s\leq\eps,\int^t_0(|\beta^1_s|^2+|\beta^2_s|^2)\dif s>18R^2\eps^{\frac{1}{4}}\right\}\\
&\subset\left\{\int^t_0|Y^1_s|^2\dif s\leq\eps,\int^t_0|\beta^1_s|^2>9R^2\eps^{\frac{1}{4}}\right\}
\cup\left\{\int^t_0|Y^2_s|^2\dif s\leq\eps,\int^t_0|\beta^2_s|^2>9R^2\eps^{\frac{1}{4}}\right\}.
\end{align*}
By Lemma \ref{Lemma2}, we have for some $C_2\in(0,1)$,
\begin{align}
\mP\left\{\int^t_0(|Y^1_s|^2+|Y^2_s|^2)\dif s\leq\eps,\int^t_0(|\beta^1_s|^2+|\beta^2_s|^2)\dif s>18R^2\eps^{\frac{1}{4}}\right\}\leq 8\exp\left\{-C_2\eps^{-\frac{1}{12}}\right\}.\label{EW11}
\end{align}
On the other hand, noticing that
$$
|uK_t|\geq 1-\int^t_0|uK_s|\cdot |\nabla b(X_s)|\dif s\stackrel{(\ref{TY1})}{\geq} 1- t\|\nabla b\|_\infty\e^{\|\nabla b\|_\infty t}\geq\tfrac{1}{2},
$$
provided $t<1\wedge(2\|\nabla b\|_\infty\e^{\|\nabla b\|_\infty})^{-1}$, we have
\begin{align}
&\mP\left\{\int^t_0(|Y^1_s|^2+|Y^2_s|^2)\dif s\leq\eps,\int^t_0(|\beta^1_s|^2+|\beta^2_s|^2)\dif s\leq18R^2\eps^{\frac{1}{4}}\right\}\no\\
&\quad\leq\mP\left\{\int^t_0(|Y^1_s|^2+|Y^2_s|^2+|\beta^1_s|^2+|\beta^2_s|^2)\dif s\leq\eps+18R^2\eps^{\frac{1}{4}}\right\}\no\\
&\quad\stackrel{(\ref{Con})}{\leq}\mP\left\{c_2\int^t_0|uK_s|^2\dif s\leq\eps+18R^2\eps^{\frac{1}{4}}\right\}\leq\mP\left\{\frac{c_2t}{4}\leq(1+18R^2)\eps^{\frac{1}{4}}\right\},\label{EW2}
\end{align}
which equals to zero provided $\eps<(\frac{c_2t}{4(1+18R^2)})^4$. If we choose
$$
C_1:=\left(\frac{c_2}{4(1+18R^2)}\right)^4,\ \ t_0:=1\wedge(2\|\nabla b\|_\infty\e^{\|\nabla b\|_\infty})^{-1},
$$
then combining (\ref{EW11}) and (\ref{EW2}), we obtain (\ref{EW3}).
\end{proof}
\bl\label{Le55}
Under {\bf (H$^\alpha_1$)} and (\ref{Con}), there exists constants $C_1,C_2\in(0,1), C_3>1$  
independent of the starting point $x$ and $t_0\in(0,1)$ such that for all $t\in(0,t_0)$ and $\eps\in(0, C_1t^{24/\alpha})$,
\begin{align}
\sup_{|u|=1}\mP(u\Sigma_t u^*\leq\eps)\leq C_3\exp\left\{-C_2\eps^{-\frac{\alpha}{72}}\right\},
\end{align}
where $\Sigma_t$ is defined by (\ref{TY7}).
\el
\begin{proof}
Noticing that
$$
u\Sigma_t u^*:=\int^t_0|uK_sA_1|^2\dif s+\int^t_0\!\!\!\int_{\Gamma_0}|uK_s A_2|^2\zeta(z) N(\dif s,\dif z),
$$
we have
\begin{align*}
\mP(u\Sigma_t u^*\leq\eps)&\leq\mP\left(u\Sigma_t u^*\leq\eps; \int^t_0(|uK_sA_1|^2+|uK_sA_2|^2)\dif s>\eps^{\frac{\alpha}{6}}\right)\\
&\quad+\mP\left(\int^t_0(|uK_sA_1|^2+|uK_sA_2|^2)\dif s\leq\eps^{\frac{\alpha}{6}}\right)\\
&\leq\mP\left(u\Sigma_t u^*\leq\eps; \int^t_0|uK_sA_1|^2\dif s>\frac{\eps^{\frac{\alpha}{6}}}{2}\right)\\
&\quad+\mP\left(u\Sigma_t u^*\leq\eps; \int^t_0|uK_sA_2|^2\dif s>\frac{\eps^{\frac{\alpha}{6}}}{2}\right)\\
&\quad+\mP\left(\int^t_0(|uK_sA_1|^2+|uK_sA_2|^2)\dif s\leq\eps^{\frac{\alpha}{6}}\right),
\end{align*}
which gives the desired estimate by Lemmas \ref{Le52} and \ref{Le53}.
\end{proof}

Now we are in a position to give:

\begin{proof}[Proof of Theorem \ref{Main2}]
By Lemma \ref{Le55} and a standard compact argument (cf. \cite[p.133 Lemma 2.31]{Nu} or \cite{Zh2}), for any $p\geq 1$, there exist constant $C_p>0$ and $\gamma(p)>0$ such that for all $t\in(0,1)$,
\begin{align}
\sup_{x\in\mR^d}\mE(\det\Sigma_t(x))^{-p}\leq C_pt^{-\gamma(p)}.\label{Es1}
\end{align}
Now, by the chain rule, we have
\begin{align*}
&\nabla^k\mE\Big((\nabla^nf)(X_t(x))\Big)
=\sum_{j=1}^k\mE\Big((\nabla^{n+j}f)(X_t(x))G_j(\nabla X_t(x),\cdots,\nabla^k X_t(x))\Big),
\end{align*}
where $\{G_j, j=1,\cdots,k\}$ are real polynomial functions.
Using (\ref{WQ2}) and as in (\ref{WQ3}) (cf. \cite[p.100, Proposition 2.1.4]{Nu}), by Lemma \ref{Le32} and H\"older's inequality, 
there exist $p_1,p_2>1$, $C>0$ independent of $x$ such that for all $t\in(0,1)$,
$$
|\nabla^k\mE((\nabla^nf)(X_t(x)))|\leq C\|f\|_\infty(\mE(\det\Sigma_t(x))^{-p_1})^{1/p_2}\leq Ct^{-\gamma(p_1)/p_2}.
$$
The proof is complete.
\end{proof}

{\bf Acknowledgements:}

The authors would like to thank Professors Zhao Dong, Feng-Yu Wang and Xuhui Peng, Lihu Xu for their interests and stimulating discussions.
Special thanks go to Lihu Xu for sending us his preprint paper \cite{Xu}.
X. Zhang is supported by NNSFs of China (Nos. 11271294, 11325105). 


\end{document}